\input amstex
\documentstyle{amsppt}

\def\id{\text{id}}
\def\qed{\qquad \vrule height6pt width5pt depth0pt}
\def\z{{\bold Z}}
\def\t{{\bold T}}
\def\oh{{\Cal O}}
\def\a{{\Cal A}}
\def\e{{\Cal E}}
\def\f{{\Cal F}}
\def\u{{\Cal U}}
\def\calt{{\Cal T}}
\def\b{{\Cal B}}
\def\calr{{\Cal R}}
\def\s{{\Cal S}}
\def\alhat{{\widehat\alpha}}
\def\Chi{{\raise2pt\hbox{$\chi$}}}
\def\restrictedto#1{\big|\lower3pt\hbox{$\scriptstyle #1$}}
\def\cstar{$C^*$}
\def\inv{^{-1}}

\input epsf         
\epsfverbosetrue    
\def\Xy{\leavevmode
 \hbox{\kern-.1em X\kern-.3em\lower.4ex\hbox{Y\kern-.15em}}}

\def\thmj{Definition 2.1}   
\def\thmmm{Definition 2.2}  
\def\thmk{Definition 2.3}   
\def\thml{Definition 2.4}   
\def\thmm{Lemma 2.5}        
\def\thmoo{Definition 2.6}      
\def\thmp{Lemma 2.7}      
\def\thmn{Definition 2.8}   
\def\thmbbb{Notes 2.9}  
\def\thmna{Definition 2.10}   
\def\thmnb{Lemma 2.11}     
\def\thmnc{Lemma 2.12}     
\def\thmo{Lemma 2.13}        
\def\thmq{Definition 2.14}  
\def\thmr{Lemma 2.15}       
\def\thms{Definitions 2.16} 
\def\thmu{Lemma 2.17}       
\def\thmv{Corollary 2.18}   
\def\thmw{Definition 2.19}   
\def\thmii{Lemma 3.1}       
\def\thmll{Definition 3.2}      
\def\thmjj{Definition 3.3}  
\def\thmpp{Lemma 3.4}   
\def\thmqq{Corollary 3.5} 
\def\thmrr{Corollary 3.6} 
\def\thmnn{Definitions 3.7} 
\def\thmzz{Proposition 3.8}    
\def\thmxx{Corollary 3.9}  
\def\thmgg{Lemma 3.10}       
\def\thmz{Corollary 3.11}    
\def\thmhh{Lemma 3.12}       
\def\thmss{Definition 3.13}  
\def\thmtt{Remarks 3.14}  
\def\thmuu{Definition 3.15} 
\def\thmvv{Definition 3.16} 
\def\thmww{Proposition 3.17} 
\def\thmx{Theorem 3.18}      
\def\thmy{Corollary 3.19}    
\def\thmaaa{Proposition 3.20}  
\def\thmaa{Lemma 4.1}       
\def\thmbb{Lemma 4.2}      
\def\thmyy{Lemma 4.3}    
\def\thmcc{Lemma 4.4}      
\def\thmdd{Lemma 4.5}      
\def\thmee{Lemma 4.6}      
\def\thmff{Theorem 4.7}    
\def\thmkk{Theorem 4.8}    

\def\displayl{2.11.1}  
\def\displaym{2.11.2}  
\def\displayo{2.11.3}  
\def\displayn{2.11.4}  
\def\displaya{3.17.1}        
\def\displayb{3.17.2}        
\def\displayc{3.17.3}        
\def\displayd{4.3.1}     
\def\displaye{4.3.2}     
\def\displayf{4.3.3}     
\def\displayg{4.3.4}     
\def\displayh{4.3.5}     
\def\displayi{4.3.6}     
\def\displayj{3.12.1}     
\def\displayk{3.12.2}     

\def\bib#1{[#1]}
\def\anantharaman{1}
\def\bkp{2}
\def\cuntzon{3}
\def\cuntzoa{4}
\def\cuntzkrieger{5}
\def\kirchberg{6}
\def\kumjianpask{7}
\def\kprr{8}
\def\lacaspielberg{9}
\def\phillips{10}
\def\renault{11}
\def\rordam{12}
\def\spielberga{13}        
\def\spielbergb{14}        
\def\spielbergc{15}        
\def\spielbergd{16}        
\def\szymanski{17}          
\def\zhang{18}

\topmatter
\title Graph-Based Models for Kirchberg Algebras
\endtitle
\author Jack Spielberg \endauthor
\address Department of Mathematics and Statistics,
Arizona State University,
Tempe, AZ  85287-1804
\endaddress
\email jack.spielberg\@asu.edu\endemail
\abstract
We give a construction of Kirchberg algebras from graphs.  
By using product graphs in the construction we are able to provide 
models for general (UCT) Kirchberg algebras while maintaining the 
explicit generators and relations of the underlying graphs.
\endabstract
\keywords simple purely infinite \cstar-algebra, K-theory, 
automorphism group, 
\break graph algebra
\endkeywords
\subjclass Primary 46L05, 46L80.  Secondary 22A22
\endsubjclass
\endtopmatter
\document

\head 1.  Introduction \endhead

The first examples of what are now termed Kirchberg algebras were
introduced in Cuntz's 1977 paper \bib\cuntzon.
(Following \bib\rordam\ we use
the term {\it Kirchberg algebra\/} for a separable nuclear simple
purely infinite \cstar-algebra.)
The brilliance of
those examples led to huge advances in the field of \cstar-algebras,
as well as deep connections with other areas of mathematics.  In a
series of papers Cuntz isolated and studied the key property of {\it
pure infiniteness\/} and its ramifications for $K$-theory (we refer to
\bib\rordam\ for a detailed bibliography).  This study was
carried further by many mathematicians, notably R\o rdam and Kirchberg.
The culmination was the classification theorem of Kirchberg
(\bib\kirchberg), also proved independently by Phillips (\bib\phillips):  
the Kirchberg algebras satisfying the universal coefficient theorem
are classified by $K$-theory.

Kirchberg algebras arise in many diffferent contexts.  As a result of 
the classification theorem, examples from different situations may be
identified by computing $K$-theory (see \bib\lacaspielberg\ for an example 
involving dynamical systems).  Alternatively, to prove a theorem about
Kirchberg algebras in general, one can choose a suitable realization
that lends itself to the problem at hand.  One of the most useful of
these has been the \cstar-algebras defined by directed graphs.  This
idea was implicitly present in Cuntz's original paper, and for finite 
irreducible graphs was explicit in the papers \bib\cuntzkrieger,
\bib\cuntzoa.  The development for arbitrary directed graphs began
with the article \bib\kprr, and has grown to be a mini-industry in
itself (the forthcoming CBMS lectures notes of Raeburn will give a
comprehensive survey).  

There have been few applications of graph
alebras to the study of Kirchberg algebras.  We mention Szyma\'nski's
proof (\bib\szymanski) that any Kirchberg algebra with 
free $K_1$-group can be realized
as the \cstar-algebra of an irreducible row-finite graph, and the
proof in \bib\spielbergb\ that Kirchberg algebras having finitely
generated $K$-theory and free $K_1$-group are semiprojective (in the
sense of Blackadar).  Graph algebras are particularly well-suited for 
such arguments in that they are defined by generators and simple,
highly flexible relations.  Their defect is clear in the
above-mentioned theorems --- graph \cstar-algebras have free
$K_1$-groups, and hence cannot be used to model general Kirchberg
algebras.  More recently, Kumjian and Pask have introduced a notion of
higher-rank graphs, or
{\it $k$-graphs\/}, and their \cstar-algebras.  These have many features in
common with ordinary (or 1-) graphs, and allow for more general 
$K_1$-groups.  However they are much less flexible than 1-graphs, and 
the theory has not yet been developed as extensively.  In particular,
the case of higher rank graphs that are not row-finite, which is
crucial for our methods, has not yet been treated.

In this paper we present examples of hybrid objects mixing elements of
$k$-graphs of different ranks.  While a general treatment seems beyond 
current technique (and perhaps not worth the considerable effort), the
special situation developed here allows us to model arbitrary
Kirchberg algebras with the same flexibility exhibited by ordinary
graph algebras.  In particular, we use this construction in
\bib\spielbergc\ and \bib\spielbergd\ to prove interesting properties 
of general (UCT) Kirchberg algebras:
\roster
\item any prime-order automorphism of the $K$-theory of a UCT
Kirchberg algebra is induced from an automorphism of the algebra
having the same order (generalizing work of
\bib\bkp);
\item a UCT Kirchberg algebra is weakly semiprojective if and only if
its $K$-groups are direct sums of cyclic groups.
\endroster
We hope that the interest of these applications will justify the work 
required to establish these models.

The rest of the paper is organized as follows.  We conclude the
introduction by recalling the basic notions of graph \cstar-algebras
in a form convenient for our purposes.  In part two we construct the 
hybrid object underlying our algebras, and use it to define an
$r$-discrete groupoid whose unit space is an appropriate set of paths.
 We remark that because the underlying object is not ``row-finite,''
 the paths we use may be finite, infinite, or semi-infinite.  In part 
 three we give generators and relations for the \cstar-algebra of this
 groupoid, prove the gauge-invariant uniqueness theorem, and show that
 the \cstar-algebra we have constructed is a UCT Kirchberg algebra. 
 Finally in part four we compute the $K$-theory of the algebra, showing
 that it is equal to the direct sum of the $K$-theory of tensor
 products of the ordinary graph algebras used in the construction of the
 underlying object.  It is here that the flexibility inherent in graph
 algebras may be used to construct our models of
 Kirchberg algebras with arbitrary
 $K$-theory.
 
 We now briefly recall the main facts 
 about (ordinary) graph algebras (see \bib\spielberga).
 A {\it directed graph\/} $E$ consists of two sets, $E^0$ (the {\it 
 vertices\/}) and $E^1$ (the {\it edges\/}), together with two maps 
 $o$, $t:E^1\to E^0$ ({\it origin\/} and {\it terminus\/}).  A {\it 
 path\/} of length $n$ in $E$ is a string $e_1e_2\cdots e_n$ of edges 
 with $t(e_i)=o(e_{i+1})$, $1\le i<n$.  We let $E^n$ denote the set of 
 paths of length $n$, and $E^*$ the set of all finite paths; 
 the origin and terminus maps extend to $E^*$ in the 
 obvious way.  For a vertex $a\in E^0$ we use the notation 
 $E^n(a)$, respectively $E^*(a)$, for the set of paths in 
 $E$ of length $n$, respectively of arbitrary length, with origin $a$.
 We let $\oh(E)$ denote the \cstar-algebra of $E$.  It 
 is the universal \cstar-algebra defined by generators $\bigl\{P_a 
 \bigm| a\in E^0\bigr\}$ and $\bigl\{S_e \bigm| e\in E^1\bigr\}$ with 
 the {\it Cuntz-Krieger relations\/}:
\roster
 \item"$\bullet$" $\bigl\{P_a \bigm| a\in E^0\bigr\}$ are pairwise 
 orthogonal projections.
 \item"$\bullet$" $S_e^*S_e=P_{t(e)}$, for $e\in E^1$.
 \item"$\bullet$" $o(e)=o(f) \Longrightarrow S_eS_e^* + S_fS_f^* \le 
 P_{o(e)}$, for $e$, $f\in E^1$ with $e\not=f$.
 \item"$\bullet$" $0 < \#\ E^1(a)<\infty \Longrightarrow P_a= 
 \sum\bigl\{S_eS_e^*\bigm| o(e)=a\bigr\}$, for $a\in E^0$.
\endroster
 \par\noindent
 (These are a variant of the relations given in \bib\spielberga, 
 Theorem 2.21.)

 The relationship between the \cstar-algebras of a graph and a 
 subgraph are crucial to our methods.  We refer to 
 \bib\spielberga.  The results are as follows.  
 Let $E$ be a graph and let $F$ be a subgraph of $E$.  We let 
 $S=S(F)$ be 
 the set of vertices in $F^0$ that do not emit more edges in $E$ than 
 in $F$.  We let $\calt\oh(F,S)$ denote the relative Toeplitz 
 Cuntz-Krieger algebra of $F$ in $E$.  It is the universal 
 \cstar-algebra defined by generators $\bigl\{P_a 
 \bigm| a\in F^0\bigr\}$ and $\bigl\{S_e \bigm| e\in F^1\bigr\}$ with 
 the relations (as above) for $\oh(F)$, modified by requiring the 
 fourth relation only if $a\in S$. Then $\calt\oh(F,S)$ is the 
 \cstar-subalgebra of $\oh(E)$ generated by the projections and partial 
 isometries associated to the vertices and edges of $F$ 
 (\bib\spielberga, Theorem 2.35).
 \smallskip
 The figures in this paper were prepared with \Xy-pic.

\head 2.  Groupoid Models for  Kirchberg algebras \endhead

We will construct models for  Kirchberg algebras
by using a mixture of directed 1-graphs and {\it 
$k$-graphs\/}, as studied by Kumjian and Pask in \bib\kumjianpask.
Since neither
the results of 
\bib\kumjianpask\ nor of \bib\spielberga\  directly apply 
in this situation, we will carry out the necessary constructions in 
detail.  Our models will consist of a sequence of
product $k$-graphs connected to 
each other by (ordinary) 1-graphs.  In fact, all details of the
argument are already present in the case of two product graphs of rank
2, and it is this case that we will treat.
The argument in the general case is essentially identical to the one we will 
give.  (We remark that in the general case one may attach other 
product $k$-graphs to $u_0$ in the same way that $E_1\times F_1$ is 
attached in what follows.)

\definition{\thmj} Let $D$ be the graph with
$$\align 
D^0&=\{u_0,\;u_1,\;a_0,\;a_1\} \\
D^1&=\{\alpha_0,\;\alpha_1\;\beta_0,\;\ldots,\;\epsilon_1\}
\endalign$$
(see Figure 1).
\enddefinition

\bigskip
\epsfysize=2.5in                   
\centerline{\epsfbox{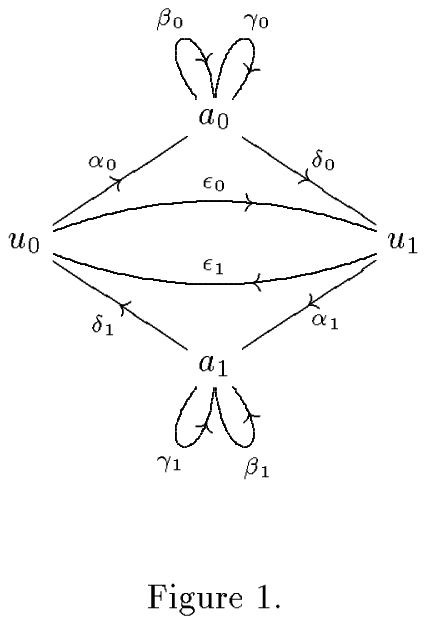}}
\bigskip

For $i=0$, 1 let $E_i$ and $F_i$ be irreducible graphs.  Let $L_i$, 
respectively $M_i$, denote the set of vertices in $E_i^0$, 
respectively $F_i^0$, emitting infinitely many edges.  We assume 
that $L_i$ and $M_i$ are nonempty.  We choose distinguished elements 
$v_i\in L_i$ and $w_i\in M_i$.   We attach the 2-graphs $E_i\times 
F_i$ to the graph $D$ by identifying $u_i\in D^0$ with $(v_i,w_i)\in 
E_i^0\times F_i^0$.  The entire object comprising $D$, $E_0\times 
F_0$, and $E_1\times F_1$ will be denoted $\Omega$.

\definition{\thmmm}
By a {\it vertex\/} we mean an element of 
$$\cup_i (E_i^0\times F_i^0) \cup  
D^0,$$ 
where we identify $u_i$ and $(v_i,w_i)$.  By an {\it edge\/} we mean 
an element of $$\bigl(\cup_i(E_i^1\times F_i^0) \cup (E_i^0 
\times F_i^1)\bigr) \cup D^1.$$
\enddefinition

\definition{\thmk} A {\it finite path element of type $D$\/} is a 
finite directed path in $D$ of non-zero length.  An {\it infinite 
path element of type $D$\/} is either an infinite directed path in 
$D$, or a finite path element of type $D$ which terminates at 
$u_0$ or $u_1$.  A {\it finite path element of type 
$(E_i,F_i)$\/} is an ordered pair $(p,q)\in E_i^*\times F_i^*$
such that $p$, $q$ are not both of length zero.
An {\it infinite path element of type $(E_i,F_i)$\/} is an 
ordered pair $(p,q)$, where $p$, 
respectively $q$, is either an infinite path or a finite path 
terminating in $L_i$, respectively in $M_i$, in $E_i$, 
respectively $F_i$,
and $p$, $q$ are not both of length zero.

For a path element $(p,q)$ of type $(E_i,F_i)$ we define {\it 
origin\/} and {\it terminus\/} by
$$\align
o(p,q)&=\bigl(o(p),o(q)\bigr) \\
t(p,q)&=\bigl(t(p),t(q)\bigr).
\endalign$$
If $(p,q)$ and $(p',q')$ are path elements of type $(E_i,F_i)$ we say 
that $(p,q)$ {\it extends\/} $(p',q')$ if  $p$ extends $p'$ and $q$ 
extends $q'$ in the usual sense of paths in a directed graph.  

A {\it finite path\/} is either a vertex, or a finite string 
$\mu_1\ldots \mu_k$ of finite path elements such that
\smallskip
\item{(i)} $t(\mu_i)=o(\mu_{i+1})$, and
\item{(ii)} $\mu_i$ and $\mu_{i+1}$ are of different types.
\smallskip
An {\it infinite path\/} is either a vertex in $L_i\times M_i$,
an infinite string of finite path 
elements satisfying the conditions (i) and (ii) above,
or a finite sequence $\mu_1\ldots \mu_{k+1}$ such that
\smallskip
\item{$\bullet$} $\mu_1\ldots \mu_k$ is a finite path
\item{$\bullet$} $\mu_{k+1}$ is an infinite path element
\item{$\bullet$} (i) and (ii) above hold.
\smallskip
\noindent
We let $X$ denote the set of all infinite paths.
We will 
use the notation $\mu\preceq \nu$ to indicate that the path 
$\nu$ extends the path  $\mu$.
\enddefinition

We wish to define a topology on $X$ making it a locally compact 
metrizable space.  First we will define the length function on paths.

\definition{\thml}  We define $\ell:\{\hbox{finite\ paths}\}\to\z_+^2$ by
$$\ell(\mu)=
 \cases
(0,0),&\text{if $\mu$ is a vertex} \\
 \bigl(\ell(p),\ell(p)\bigr),&\text{if $\mu=p$ 
 is a finite path element of type $D$} \\
\bigl(\ell(p),\ell(q)\bigr),&\text{if $\mu=(p,q)$ 
 is a finite path element of type $(E_i,F_i)$.} \\
 \endcases$$
 (In the right hand side above we have used the symbol $\ell$ also for 
 the usual length function on paths in a 1-graph.  There should be no 
 confusion resulting from this usage.)
If $\mu=\mu_1\ldots \mu_k$ is a finite path, we define
$$\ell(\mu)=\sum_{i=1}^k\ell(\mu_k).$$
We also will use the extension of $\ell$ to  
infinite paths, with values in $\bigl(\z_+\cup\{\infty\}\bigr)^2$, defined in the 
obvious way.
\enddefinition

\proclaim{\thmm} $\ell$ is additive with respect to concatenation of paths.
\endproclaim

\demo{Proof} The proof is left to the reader. \qed
\enddemo

\definition\thmoo  Let $\mu=\mu_1\mu_2\cdots\mu_r$ be a finite path,
decomposed as a string of finite path elements satisfying \thmk(i) and
(ii).  We write $|\mu|=r$, the number of finite path elements in $\mu$.
\enddefinition

We wish to define a groupoid having having $X$
as its unit space.  
Since $\Omega$ is not a higher rank graph, 
the space of paths does not have the 
factorization property (Definition 1.1 of \bib\kumjianpask).
However we do have the following simple observation.

\proclaim{\thmp} Let $\mu$ and $\mu'$ be finite paths, and let 
$x$, $x'$ be infinite paths, with $\mu x=\mu' x'$.  Let 
$m=\ell(\mu)\vee \ell(\mu')$. (Here $\vee$ represents the usual 
lattice join in $\z^2$.)
Then there is a factorization $\mu x=\nu y$ with $\ell(\nu)=m$.
\endproclaim

\demo{Proof}  The proof is accomplished by an easy
induction on $n=\max\{|\mu|,\,|\mu'|\}$.  \qed
\enddemo

\definition{\thmn} Let $\mu=\mu_1\ldots \mu_k$ be a finite path
decomposed as a string of finite path elements satisfying \thmk(i) and
(ii).  Let 
$$\align
Z(\mu)&=\bigl\{\sigma_1\sigma_2\cdots\in X : \sigma_i=\mu_i \text{ for } i<k,\ 
\mu_k \preceq \sigma_k\bigr\}. \\
\intertext{We define certain subsets of $Z(\mu)$.}
V(\mu)&=
\cases
   Z(\mu),&\text{if $t(\mu)\not\in\{u_0,u_1\}$} \\
   Z(\mu)\setminus Z(\mu \alpha_i)\setminus Z(\mu\epsilon_i),
 &\text{ if } t(\mu)=u_i. \\
\endcases \\
\intertext{If $t(\mu)=(y,z)\in E^0_i\times F_i^0$, let
$B\subseteq E^1(y)$ and $C\subseteq E^1(z)$ be finite subsets 
such that if $y\not\in L_i$ then $B=\emptyset$, and if
$z\not\in M_i$ then $C=\emptyset$.  Then we define:}
V(\mu;B,C)&=
  V(\mu)\setminus \bigcup_{e\in B} Z\bigl(\mu(e,z)\bigr) 
  \setminus \bigcup_{f\in C} Z\bigl(\mu(y,f)\bigr)
\endalign$$
\enddefinition

Note that for $t(\mu)\in E_i^0\times F_i^0$, 
$V(\mu)=V(\mu;\emptyset,\emptyset)$.
Also note that if $t(\mu)=u_i$ for some $i$, then every element of
$V(\mu;B,C)\setminus\{\mu\}$ extends $\mu$ farther into $E_i\times 
F_i$.

\remark{Remarks} We mention that the topology we intend to define on 
$X$ requires that we be able to `block' finitely many edges 
originating at a vertex of $\Omega$ that emits infinitely many edges. 
This explains our use of the notation $V(\mu;B,C)$.  At a vertex 
emitting only finitely many edges, we can handle the situation by 
considering explicitly the finitely many possible extensions by an 
edge.  Since we do not want to block all edges originating at such a 
vertex, we do not allow $B$ and/or $C$ to be nonempty at such a 
vertex.
\endremark

\example{\thmbbb} We note the following disjoint unions of sets.
\roster
\item If $t(\mu)=u_i$ for some $i$, then for $B$ and $C$ nonempty,
$$\align
Z(\mu)&=V(\mu)\cup V(\mu\alpha_i) \cup Z(\mu\epsilon_i).
\endalign$$
\item If $t(\mu)=(y,z)\in E_i^0\times F_i^0$, then
$$\align
V(\mu)
&= V(\mu;B,C) \\
&\qquad \cup \bigcup\Bigl\{Z\bigl(\mu(e,f)\bigr)
   : e\in B,\;f\in C\Bigr\} \\
\allowdisplaybreak
&\qquad \cup \bigcup_{e\in B} V\bigl(\mu(e,z);\emptyset,C\bigr) \\
&\qquad\quad \cup \bigcup\Bigl\{ V\bigl(\mu(e,z)\alpha_i\bigr) \cup
   Z\bigl(\mu(e,z)\epsilon_i\bigr) : e\in B,\; t(e,z)=u_i\Bigr\} \\
\allowdisplaybreak
&\qquad \cup \bigcup_{f\in C} V\bigl(\mu(y,f);B,\emptyset\bigr) \\
&\qquad\quad \cup \bigcup\Bigl\{ V\bigl(\mu(y,f)\alpha_i\bigr) \cup
   Z\bigl(\mu(y,f)\epsilon_i\bigr) : f\in C,\; t(y,f)=u_i\Bigr\}, \\
V(\mu)
&=V(\mu;B,\emptyset)  \cup \bigcup_{e\in B} Z\bigl(\mu(e,z)\bigr), \\
V(\mu)
&=V(\mu;\emptyset,C)  \cup \bigcup_{f\in C} Z\bigl(\mu(y,f)\bigr).
\endalign$$
\endroster
\endexample

\definition{\thmna} We let $\b$ denote the collection of all 
$V(\mu;B,C)$ and $Z(\mu)$.  
We let $\a$ denote the collection of all finite 
disjoint unions of sets in $\b$ (hence also $\emptyset\in\a$.)
\enddefinition

\proclaim{\thmnb} The collections $\a$ and $\b$ have
the following properties:
\roster
\item $\emptyset\not\in\b$.
\item $\b$ is countable.
\item $\a$ is a ring of sets.
\item The intersection of any decreasing sequence of sets in 
$\b$ is nonempty.
\endroster
\endproclaim

\demo{Proof} A bit of thought will likely convince the reader 
that this is elementary.  However we have found the details slightly
messy, and we warn that the following is a bit of a slog.
(1) and (2) are clear.  For (3), we must show that if $A$, $A'\in\b$
then $A\cap A'$, $A\setminus A'\in\a$.  We will first deal with
intersections.
Let $A$, $A'\in\b$.  First suppose
that $A=Z(\mu)$ and $A'=Z(\mu')$.  If $\mu\preceq\mu'$ then $A\cap
A'=A'$.  So suppose that $\mu$, $\mu'$ are not comparable.  
We then have
$$\aligned
\mu&=\mu_1\cdots\mu_{k-1}\mu_k\sigma \\
\mu'&=\mu_1\cdots\mu_{k-1}\mu_k'\sigma',
\endaligned \tag\displayl$$
decomposed as a string of finite path elements satisfying \thmk(i) and
(ii),
where $\mu_k$ and $\mu_k'$ are not comparable.  We have $Z(\mu)\cap
Z(\mu')=\emptyset$, unless $\sigma$ and $\sigma'$ are absent, $\mu_k$ 
and $\mu_k'$ are of type $(E,F)$, and
$$\aligned
\mu_k&=\widetilde{\mu}(e_1\cdots e_m,z) \\
\mu_k'&= \widetilde{\mu}(y,f_1\cdots f_n)
\endaligned \tag\displaym$$
with $m$ and $n$ both nonzero.  In this case $Z(\mu)\cap Z(\mu') =
Z(\mu'')$, where 
$$\mu''=\mu_1\cdots\mu_{k-1}\widetilde{\mu} (e_1\cdots
e_m,\,f_1\cdots f_n).\tag\displayo$$

Next suppose that $A=V(\mu;B,C)$ and $A'=Z(\mu')$.  
First consider the situation where $\mu$ and $\mu'$ are comparable.  If
$\mu'\preceq\mu$, then we have $A\subseteq Z(\mu)\subseteq A'$.
If $\mu\precneqq\mu'$, write $\mu$ and $\mu'$ as
$$\aligned
\mu&=\mu_1\cdots\mu_k \\
\mu'&=\mu_1\cdots\mu_{k-1}\mu_k'\sigma_1\cdots\sigma_j, \\
\mu_k'&=\mu_k\widetilde{\mu}.
\endaligned \tag\displayn$$
There are two cases.

\remark{Case (i)}
Suppose that $\mu_k=\mu_k'$.  Then we must have $j>0$.  If $\sigma_1$ 
is of type $D$ then $V(\mu;B,C)\subseteq V(\mu)$, while
$Z(\mu')\subseteq Z(\mu\alpha_i)\cup Z(\mu\epsilon_i)$ for some $i$.
Hence $A\cap A'=\emptyset$ (\thmbbb(1)).  
If $\sigma_1$ is of type $(E,F)$ then we have
$$A\cap A'=
\cases
Z(\mu'),&\text{if $\sigma_1$ obeys the restrictions imposed by $B$ and 
$C$} \\
\emptyset,&\text{otherwise.}
\endcases$$
\endremark

\remark{Case (ii)}
Suppose that $\mu_k\precneqq\mu_k'$.  If $\mu_k$ and $\mu_k'$ are of type 
$D$ then we have
$$A\cap A'=
\cases
A',&\text{if $t(\mu)\in\{a_0,\,a_1\}$} \\
\emptyset,&\text{otherwise,}
\endcases$$
for the same reason as in case (i).
If $\mu_k$ and $\mu_k'$ are of type $(E,F)$ we assume that
$\widetilde{\mu}$ obeys the restrictions imposed by $B$ and $C$ (since
otherwise we would have $A\cap A'=\emptyset$).  If $j>0$ or if
$\ell(\widetilde{\mu})\ge(1,1)$ then $A'\subseteq A$.  Otherwise we
have $j=0$ and $\ell(\widetilde{\mu})=(m,n)$ with exactly
one of $m$, $n$
equal to zero.  For definiteness suppose $n=0$.  Then 
$$A\cap A'=
\cases
V(\mu';\emptyset,C),&\text{if $t(\mu')\not\in\{u_0,\,u_1\}$} \\
V(\mu';\emptyset,C)\cup V(\mu'\alpha_i)\cup Z(\mu'\epsilon_i),&
\text{if $t(\mu')=u_i$ for some $i$}.
\endcases$$
\endremark

Now consider the situation where $\mu$ and $\mu'$ are not comparable. 
Write $\mu$ and $\mu'$ as in (\displayl).  Since $A\subseteq Z(\mu)$, 
we have $A\cap A'=\emptyset$ unless we are in the situation of
(\displaym).  In this case let $\mu''$ be as in (\displayo).  We
have
$$\align
A\cap A'
&=V(\mu;B,C)\cap Z(\mu') \\
&=V(\mu;B,C)\cap Z(\mu'') \\
&=
\cases
\emptyset,&\text{if $f_1\in C$} \\
V(\mu'';B,\emptyset),&\text{if $f_1\not\in C$, $t(\mu'')\not\in
\{u_0,\,u_1\}$} \\
V(\mu'';B,\emptyset) \cup V(\mu''\alpha_i) \cup Z(\mu''\epsilon_i),&
\text{if $f\not\in C$, $t(\mu'')=u_i$, some $i$.}
\endcases
\endalign$$

Finally we let $A=V(\mu;B,C)$ and $A'=V(\mu';B',C')$.  Again, we first
consider the situation where $\mu$ and $\mu'$ are comparable.  If
$\mu=\mu'$ then we have
$$A\cap A'=V(\mu;B\cup B',C\cup C').$$
Suppose now that $\mu\precneqq\mu'$  Write $\mu$ and $\mu'$ as in
(\displayn).  There are two cases.

\remark{Case (iii)}
Suppose that $\mu_k=\mu_k'$.  Then $j>0$.  If $\sigma_1$ is of type
$D$ then $A\cap A'=\emptyset$, since $A'\subseteq Z(\mu')$.  If
$\sigma_1$ is of type $(E,F)$, and if
$\ell(\sigma_1,\ldots,\sigma_j)\ge(1,1)$, then
$$A\cap A'=
\cases
A',&\text{if $\sigma_1$ obeys the restrictions imposed by $B$ and $C$} \\
\emptyset,&\text{otherwise.}
\endcases$$
If $j=1$ and $\ell(\sigma_1)=(m,n)$ where exactly one of $m$, $n$ is
nonzero, suppose without loss of generality that $n=0$.  Then
$$A\cap A'=
\cases
V(\mu';B',C\cup C'),&\text{if $\sigma_1$ obeys the 
restrictions imposed by $B$ and $C$} \\
\emptyset,&\text{otherwise.}
\endcases$$
\endremark

\remark{Case (iv)} (This case is nearly identical to the earlier case
(ii).)
Suppose that $\mu_k\precneqq\mu_k'$.  If $\mu_k$ and $\mu_k'$ are of type 
$D$ then we have
$$A\cap A'=
\cases
A',&\text{if $t(\mu)\in\{a_0,\,a_1\}$} \\
\emptyset,&\text{otherwise.}
\endcases$$
If $\mu_k$ and $\mu_k'$ are of type $(E,F)$ we assume that
$\widetilde{\mu}$ obeys the restrictions imposed by $B$ and $C$ (since
otherwise we would have $A\cap A'=\emptyset$).  If $j>0$ or if
$\ell(\widetilde{\mu})\ge(1,1)$ then $A'\subseteq A$.  Otherwise we
have $j=0$ and $\ell(\widetilde{\mu})=(m,n)$ with exactly
one of $m$, $n$
equal to zero.  For definiteness suppose $n=0$.  Then 
$$A\cap A'=V(\mu';B',C\cup C').$$
\endremark

Now consider the case where $\mu$ and $\mu'$ are not comparable. 
Write $\mu$ and $\mu'$ as in (\displayl).  As before, $A\cap
A'=\emptyset$ unless we have (\displaym) with $\sigma$ and $\sigma'$ 
absent.  In this case we have
$$A\cap A'=V(\mu'',B,C'),$$
where $\mu''$ is as in (\displayo).
This finishes the proof that $\a$ is closed under intersections.

We now show that $A\setminus A'\in\a$.  We first suppose that
$A=Z(\mu)$ and $A'=Z(\mu')$.  If $\mu\preceq\mu'$ let $\mu$ and $\mu'$
be as in (\displayn).  We have a disjoint union:
$$A\setminus A'=\bigl(Z(\mu)\setminus Z(\mu\widetilde{\mu})\bigr) \cup
\bigcup_{0\le i<j}\bigl( Z(\mu\widetilde{\mu}\sigma_1\cdots\sigma_i)
\setminus Z(\mu\widetilde{\mu}\sigma_1\cdots\sigma_{i+1})\bigr).$$
Thus it suffices to consider $Z(\mu\sigma)\setminus Z(\mu)$.  If
$\sigma=e_1\cdots e_m$ is of type $D$, then
$$Z(\mu\sigma)\setminus Z(\mu) = \bigcup_{0\le i<m} \bigl( Z(\mu
e_1\cdots e_i) \setminus Z(\mu e_1\cdots e_{i+1})\bigr),$$
and for $e\in D^1$ we have
$$Z(\tau e)\setminus Z(\tau)= 
\cases
V(\tau)\cup Z(\tau e'),\quad\text{if $t(\tau)=u_i$, where $\{e'\} =
\{\alpha_i,\,\epsilon_i\}\setminus\{e\}$} \\
\bigcup\bigl\{ Z(\tau f) : f\in\{\beta_i,\,\gamma_i,\,\delta_i\}
\setminus\{e\}\bigr\},\quad\text{if $t(\tau)=a_i$.}
\endcases$$
If $\sigma=(e_1\cdots e_m,f_1\cdots f_n)$ is of type $(E,F)$, then
$$\align
Z(\mu\sigma)&\setminus Z(\mu) = 
\bigcup\Bigl\{ V(\mu\nu\alpha_i) \cup Z(\mu\nu\epsilon_i) :
o(\sigma)\precneqq\nu\precneqq\sigma,\ 
t(\nu)=u_i,\ i\in\{0,1\} \Bigr\} \\
&\cup \Bigl\{ V\bigl(\mu(e_1\cdots e_i,\,f_1\cdots f_j); \{e_{i+1}\},
\{f_{j+1}\}\bigr) : \\
&\qquad \qquad 0\le i\le m,\ 0\le j\le n,\ i+j<m+n\Bigr\},
\endalign$$
where in the second summation we let $\{e_{m+1}\}$ and $\{f_{n+1}\}$
denote the empty set.

If $\mu$ and $\mu'$ are not comparable, an earlier part of the proof
implies that $A\cap A'=\emptyset$ unless $\mu$ and $\mu'$ are as in
(\displayl) and (\displaym), with $\sigma$ and $\sigma'$ absent.   
Then
$$A\setminus A'=\bigcup_{0\le j<n} Z\bigl(\mu(y,f_1\cdots f_j)\bigr)
\setminus Z\bigl(\mu(y,f_1\cdots f_{j+1})\bigr),$$
and the differences in the above union are treated by
$$Z\bigl(\tau(y,f)\bigr)\setminus Z(\tau) =
\cases
V(\tau\alpha_i)\cup Z(\tau\epsilon_i) \cup V(\tau;\emptyset,\{f\}),
&\text{if $t(\tau)=u_i$} \\
V(\tau;\emptyset,\{f\}),&\text{if $t(\tau)\not\in\{u_0,u_1\}$}.
\endcases$$

Next suppose that $A=V(\mu;B,C)$ and $A'=Z(\mu')$. We have a disjoint 
union:
$$A'\setminus A = \bigl(Z(\mu')\setminus Z(\mu)\bigr) \cup
\Bigl(
Z(\mu') \cap\bigl(Z(\mu)\setminus V(\mu;B,C)\bigr)
\Bigr).$$
We have already proved that the first piece is in $\a$.  We saw in
\thmbbb(1) that $Z(\mu)\setminus V(\mu;B,C)\in\a$.  Thus $A'\setminus 
A\in\a$ by what we have proved for intersections of sets in $\b$.  A
similar argument applies to $A\setminus A'$ in view of
$$A\setminus A'=\bigl(Z(\mu)\setminus Z(\mu')\bigr) \cap V(\mu;B,C).$$
In the case $A=V(\mu;B,C)$ and $A'=V(\mu';B',C')$, a similar argument 
together with the above case applies:
$$A\setminus A' = \bigl(Z(\mu)\setminus V(\mu';B',C')\bigr) \cap
V(\mu;B,C).$$
This completes the proof of item (3) of the lemma.

Finally we prove (4).  Let $A_1\supseteq A_2\supseteq \cdots$, with
$A_i=V(\mu_i;B_i,C_i)$ or $Z(\mu_i)$.  From the proof of (3) we have
that $\mu_1\preceq\mu_2\preceq\cdots$.  Then
$x=\lim_i\mu_i\in\cap_1^\infty A_i$. \qed

\enddemo

\proclaim{\thmnc} Let $\e\subseteq\a$ have the finite intersection 
property.  Then $\bigcap\e\not=\emptyset$.
\endproclaim

\demo{Proof} Since $\a$ is countable, we may list the elements 
or $\e$ as $A_1$, $A_2$, $\ldots$.  Thus our assumption on $\e$ takes 
the form $\bigcap_{i=1}^p A_i\not=\emptyset$ for all $p$.
We will construct
elements $A_1'$, $A_2'$, $\ldots \in \b$ such that $A_i'\subseteq A_i$
and $A_1'\supseteq A_2'\supseteq\cdots$.  Then the lemma will follow
by \thmnb(4).

Let $A_1=\bigcup_{j=1}^{k_1}A_{1j}$ be written as a 
disjoint union of elements of $\b$.  We 
claim that there exists $j_1$ such that for all $p\ge2$,
$$A_{1,j_1}\cap \bigcap_{i=2}^p A_p \not=\emptyset.$$
For if not, then for all $j=1$, $\ldots$, $k_1$ there exists $p_j$ 
such that
$$A_{1j}\cap \bigcap_{i=2}^{p_j} A_i = \emptyset.$$
Let $p=\max\{p_1,\ldots,p_{k_1}\}$.  Then 
$A_{1j}\cap \bigcap_{i=2}^p A_i=\emptyset$ for $j=1$, $\ldots$, $k_1$.  
Hence $A_1\cap \bigcap_{i=2}^p A_i =\emptyset$, 
a contradiction.  We set $A_1'=A_{1,j_1}$.

Now suppose inductively that we have found 
$A'_1$, $A'_2$, $\ldots$, $A'_{n-1}\in\b$ such that 
$A_i'\subseteq A_i$, $A_{i-1}'\supseteq A_i'$, and 
$$A_{n-1}'\cap \bigcap_{i=n}^p A_i\not=\emptyset,\quad p\ge n.$$
For $i\ge n$ let $A_i''=A_i\cap A_{n-1}'$.  Then $A_i''\subseteq
A_{n-1}'$ and $\bigcap_{i=n}^p A_i''\not=\emptyset$ for all $p\ge n$.
Let $A_n''=\bigcup_{j=1}^{k_n} A_{nj}''$ 
be written as a disjoint union of
elements of $\b$.  We claim that there exists $j_n$ such that for all 
$p\ge n+1$, 
$$A_{n,j_n}''\cap \bigcap_{i=n+1}^p A_i''\not=\emptyset.$$
For suppose not.  Then for all $j=1$, $\ldots$, $k_n$ there exists 
$p_j\ge n+1$ such that
$$A_{nj}'' \cap\bigcap_{i=n+1}^{p_j} A_i''=\emptyset.$$
Let $p=\max_{1\le j\le k_n} \{p_j\}$.  Then
$A_{n,j}''\cap \bigcap_{i=n+1}^p A_i''=\emptyset$ for $j=1$, $\ldots$,
$k_n$.  Hence $A_n''\cap \bigcap_{i=n+1}^p A_i'' =\emptyset$,
a contradiction.  Therefore $j_n$ exists as claimed, and we set
$A_n'=A_{n,j_n}''$. \qed

\enddemo

\proclaim{\thmo} The collection, $\b$, of all $V(\mu;B,C)$ is a 
base for a locally compact metrizable topology on $X$.
\endproclaim

\demo{Proof} It follows from \thmnb (2) and (3) that $\b$ is a 
base for a second countable topology on $X$.  Thus metrizability will 
follow from local compactness and the Hausdorff property.  To 
establish local compactness we must show that if $A\in\b$ is covered 
by $\u\subseteq\b$, then $A$ is finitely covered by $\u$.  Suppose 
not.  Then for every finite subcollection $\f\subseteq\u$, 
$$A\setminus \bigcup\f = \bigcap_{U\in\f}(A\setminus 
U)\not=\emptyset.$$
Let $\e=\{A\setminus U\bigm| U\in\u\}\subseteq\a$.  Since $\e$ has 
the finite intersection property, \thmnc\ implies that 
$\bigcap\e\not= \emptyset$.  But $\bigcap\e=A\setminus \bigcup\u = 
\emptyset$, a contradiction.

To verify the Hausdorff property, let $x\not=x'$ in $X$.  If 
$x$ and $x'$ are comparable, say $x\preceq x'$.  Then there exists a
finite path element $\mu$ such that $\mu\preceq x$ and $\mu\npreceq
x'$.  Then $Z(\mu)$ and $Z\bigl(o(x)\bigr)\setminus Z(\mu)$ are
disjoint neighborhoods of $x$ and $x'$.
If $x$ and $x'$ are not comparable, write $x=\nu x_1$ 
and $x'=\nu x_1'$, where $x_1$ and $x_1'$ have no common 
nonzero initial subpath.  Let $\nu_1$ and $\nu_1'$ be nonzero initial 
finite subpaths of $x_1$ and $x_1'$.  Then $Z(\nu\nu_1)$ and
$Z(\nu\nu_1')\setminus Z(\nu\nu_1)$ are disjoint neighborhoods of $x$ 
and $x'$.  \qed
\enddemo

\definition{\thmq} Let $G$ be the set of triples $(x,n,y)$ in 
$X\times \z^2\times X$ such that there exists $z\in X$ and 
decompositions $x=\mu z$, $y=\nu z$ with $\ell(\mu)-\ell(\nu)=n$.
\enddefinition

\proclaim{\thmr} $G$ is a groupoid with the operations
$$\align
(x,n,y)(y,m,z)&=(x,n+m,z) \\
(x,n,y)\inv&=(y,-n,x)
\endalign$$
\endproclaim

\demo{Proof} It suffices to show that if $(x,n,y)$ and $(y,m,z)$ are 
in $G$ then so are $(x,n+m,z)$ and $(y,-n,x)$.  It is clear that 
$(y,-n,x)\in G$.  Let $x=\mu\zeta$, $y=\nu\zeta=\nu'\zeta'$, and 
$z=\sigma\zeta'$, where $\zeta$, $\zeta'\in X$ and $\mu$, $\nu$, $\nu'$, 
$\sigma$ are finite paths with $\ell(\mu)-\ell(\nu)=n$ and
$\ell(\nu')-\ell(\sigma)=m$.  Applying \thmp\ to the equality 
$\nu\zeta=\nu'\zeta'$ we have a decomposition $\nu\zeta=\lambda\xi$, 
where $\xi\in X$ and $\lambda$ is a finite path such that 
$\ell(\lambda) = \ell(\nu)\vee\ell(\nu')$.  Then $\lambda=\nu\eta = 
\nu'\eta'$ for some finite paths $\eta$ and $\eta'$.  We obtain
$$\align
\nu\zeta&= \lambda\xi=\nu\eta\xi,\text{ so }\zeta=\eta\xi, \\
\nu'\zeta'&= \lambda\xi=\nu'\eta'\xi, \text{ so } \zeta'=\eta'\xi.
\endalign$$
Then $x=\mu\zeta=\mu\eta\xi$, $z=\sigma\zeta'=\sigma\eta'\xi$, and
$$\align 
\ell(\mu\eta)-\ell(\sigma\eta')&= 
\ell(\mu)+\ell(\eta)-\ell(\eta')-\ell(\sigma) \\
&=\ell(\mu)+\ell(\nu')-\ell(\nu)-\ell(\sigma), \text{ since } 
\nu\eta=\nu'\eta', \\
&=n+m. \qed
\endalign$$
\enddemo

\definition{\thms} Consider the collection of
subsets of $G$ of the forms
$$\align
U(\mu_1,\mu_2)&=
\Bigl(Z(\mu_1)\times\bigl\{\ell(\mu_1)-\ell(\mu_2)\bigr\} \times
Z(\mu_2)\Bigr)\cap G, \\
U_0(\mu_1,\mu_2;B,C)&=
\Bigl(V(\mu_1;B,C)\times\bigl\{\ell(\mu_1)-\ell(\mu_2)\bigr\}\times 
V(\mu_2;B,C)\Bigr)\cap G,
\endalign$$
where $\mu_1$, $\mu_2$ are finite paths with $t(\mu_1)=t(\mu_2)$.
We will also write $U_0(\mu_1,\mu_2)$ for 
$U_0(\mu_1,\mu_2;\emptyset,\emptyset)$.
It is easy to check that this collection
is a base (of compact-open $G$-sets) for a locally compact, 
Hausdorff, totally disconnected topology making $G$ into an 
$r$-discrete groupoid.  The map
$c:(x,n,y)\in G\mapsto n\in\z^2$ is clearly a continuous cocycle.
\enddefinition

\proclaim{\thmu} $G$ is topologically free, minimal, and locally 
contractive.
\endproclaim

\demo{Proof} The proof is virtually identical to the (easy) proofs for 
graph algebras in \bib\spielberga. \qed
\enddemo

\proclaim{\thmv} $C^*_r(G)$ is simple and purely infinite.
\endproclaim

\demo{Proof} This follows from \bib\renault\ and \bib\anantharaman\ 
(see also \bib\lacaspielberg). \qed
\enddemo

We next define the gauge action of $\t^2$ on $C^*(G)$.  We will use the
notation $\zeta_m$, $m\in\z^2$, for the characters of $\t^2$:
$\zeta_m(z)=z_1^{m_1} z_2^{m_2}$.

\definition{\thmw} The gauge action $\alpha:\t^2\to 
Aut\bigl(C^*(G)\bigr)$ is dual to the 
cocycle $c$.  Thus for $f\in C_c(G)$, $\alpha_z(f)(g) = 
\zeta_{c(g)}(z)f(g)$.
\enddefinition

\head 3.  Generators and Relations \endhead

For a finite path $\mu$ we let $s_\mu\in C_c(G)$ denote the partial 
isometry
$$s_\mu=\Chi_{U(\mu,t(\mu))}.$$

\proclaim\thmii $C^*(G)$ is 
generated by the set of all $s_\mu$.  Moreover, if $t(\mu_1)=o(\mu_2)$
then $s_{\mu_1\mu_2}=s_{\mu_1}s_{\mu_2}$.
\endproclaim

\demo{Proof} We first note that
$$U(\mu_1,\mu_2)=U\bigl(\mu_1,t(\mu_1)\bigr)\cdot 
U\bigl(\mu_2,t(\mu_2)\bigr)\inv,$$
so that $\Chi_{U(\mu_1,\mu_2)}=s_{\mu_1}s_{\mu_2}^*$. From \thmbbb(1) 
we have
$$\Chi_{U_0\bigl(\mu,t(\mu)\bigr)}=
\cases
s_\mu,&\text{if $t(\mu)\not\in\{u_0,u_1\}$} \\
s_\mu - s_{\mu\alpha_i} - s_{\mu\epsilon_i},&\text{if $t(\mu)=u_i$}.
\endcases$$
It follows from \thmbbb(2) that 
$\Chi_{U_0(\mu,t(\mu);B,C)}$ is in the span of the $s_\mu$.
Finally we note that
$$\align
U_0(\mu_1,\mu_2;B,C)&=U_0\bigl(\mu_1,(y,z);B,C\bigr)\cdot
U_0\bigl(\mu_2,(y,z);B,C\bigr)\inv, \\
\intertext{and hence that}
\Chi_{U_0(\mu_1,\mu_2;B,C)}
&=\Chi_{U_0(\mu_1,(y,z);B,C)}\cdot
\bigl(\Chi_{U_0(\mu_2,(y,z);B,C)}\bigr)^*
\endalign$$
is in the span of the $s_\mu$.

The last claim follows from the fact that 
$$U\bigl(\mu_1\mu_2,t(\mu_2)\bigr)=U\bigl(\mu_1,t(\mu_1)\bigr) \cdot
U\bigl(\mu_2,t(\mu_2)\bigr). \qed$$
\enddemo

In the sequel we will let $A$ denote $C^*(G)$.
We wish to give a presentation of $A$ by generators and relations.  
For this we recall the hybrid graph structure of $\Omega$ given in 
\thmmm.  

\definition\thmll We write 
$$\Omega^{(i,j)} = \{\mu \text{ a finite path } : \ell(\mu)=(i,j)\}$$
(cf. \thmk).  Thus $\Omega^{(0,0)}$ is the set of vertices in
$\Omega$, while the set of edges is
$\Omega^{(1,0)}\cup\Omega^{(0,1)}\cup D^1$.  We write
$$\Omega^*=\bigcup_{i,j}\Omega^{(i,j)}$$
for the set of all finite paths.  For $x\in\Omega^{(0,0)}$ 
we will write $\Omega^{(i,j)}(x)=\{\mu\in\Omega^{(i,j)} : o(\mu)=x\}$.

We remark that this notation is patterned on that of \bib\kumjianpask.
Note that the ``origin'' of a path in a graph corresponds to its
``range'' when it is thought of as a morphism in a small category.
\enddefinition

\definition\thmjj We let $\s$ 
denote the  set of symbols 
$$\bigl\{P_x\bigm|\ x \text{ is a vertex}\bigr\} \cup
\bigl\{S_y\bigm|\ y \text{ is an edge}\bigr\}.$$
We let $\calr$ denote the following set of relations on $\s$:
\smallskip\noindent
\item{(i)} $P_x$ is a projection for every vertex $x$,  $S_y$ is a 
partial isometry for every edge $y$.
\item{(ii)} For every $a\in E^0_i$, the projections for $\{a\}\times 
F_i^0$ and the partial isometries for $\{a\}\times F_i^1$ satisfy the 
Cuntz-Krieger relations corresponding to the graph $F_i$ (see 
the discussion at the end of the section 1).
\item{(ii')} For every $b\in F^0_i$, the projections for 
$E_i^0\times \{b\}$ and the partial isometries for 
$E_i^1\times \{b\}$ satisfy the 
Cuntz-Krieger relations corresponding to the graph $E_i$.
\item{(iii)} The projections for $D^0$ and 
the partial isometries for $D^1$ satisfy the Toeplitz-Cuntz-Kriger
relations corresponding to the graph $D$ and the vertices
$\{a_0,a_1\}$.
\item{(iv)} If $\mu$ and $\nu$ are edges of types $D$ 
and $E_i\times F_i$, respectively, then $S_\mu^* S_\nu=0$.
\item{(v)} (2-graph structure of $E_i\times F_i$.)
For all $e\in E_i^1$ and $f\in F_i^1$ we have
$$\align
S_{(o(e),f)}\,S_{(e,t(f))} &= 
S_{(e,o(f))}\,S_{(t(e),f)} \\
S_{(t(e),f)}\,S_{(e,t(f))}^* &= 
S_{(e,o(f))}^*\, S_{(o(e),f)}.
\endalign$$
\enddefinition

We let $\Theta=C^*\langle\s,\calr\rangle$ 
denote the universal \cstar-algebra given by these 
generators and relations.  For a finite path written as a product of
edges:  $\mu=y_1y_2\cdots y_k$, 
we let $S_\mu$ denote the product $S_{y_1}S_{y_2}\cdots S_{y_k}$
(it follows from the relation 
(v) that this is unambiguous).  It is easily seen from \thmn\ and
\thms\ that $S_\mu\mapsto s_\mu$ determines a surjective 
$*$-homomorphism, $\pi$, of $\Theta$ onto $A$.  
We will show below
(\thmy) that $\pi$ is an isomorphism.  First we need to study $\Theta$
more closely.

For the next lemma recall from \thmoo\ the notation $|\mu|$ for the
number of finite path elements in $\mu$.

\proclaim{\thmpp}  Let $\mu$, $\nu\in\Omega^*$.  Suppose that $S_\mu^*
S_\nu\not=0$.  Then
\roster
\item $t(\mu)=t(\nu)$.
\item If $|\mu|>|\nu|$ then $\nu\preceq\mu$.
\item If $|\mu|<|\nu|$ then $\mu\preceq\nu$.
\item If $|\mu|=|\nu|=n$, then $\mu$ and $\nu$ are decomposed into
finite path elements satisfying \thmk\ (i) and (ii), as
$$\align
\mu&=\mu_1\mu_2\cdots\mu_n \\
\nu&=\mu_1\mu_2\cdots\mu_{n-1}\nu_n,
\endalign$$
where $\mu_n$ and $\nu_n$ are of the same type.  Moreover, if $\mu_n$
and $\nu_n$ are of type $D$ then one extends the other, while if they 
are of type $E_i\times F_i$ then in each coordinate one extends the
other.
\endroster
\endproclaim

\demo{Proof} (1) follows immediately from \thmk\ and \thmjj(iv).
We prove (2) -- (4) by 
induction on $n=\max\{|\mu|,|\nu|\}$.  The lemma is easily verified
if $n=1$.  Suppose it is true if $\max\{|\mu|,|\nu|\}<n$ where $n>1$. 
Assume that $\max\{|\mu|,|\nu|\}=n$; say $|\mu|=n\ge|\nu|$.
Write $\mu=\mu_1\cdots\mu_n$ and $\nu=\nu_1\cdots\nu_s$.  We claim
that $\nu_1\preceq\mu_1$.  We know that $S_{\mu_1}^* S_{\nu_1}\not=0$.
If $\mu_1$ and $\nu_1$ are of type $D$, by the inductive hypothesis
we must have $\mu_1\preceq\nu_1$ or $\nu_1\preceq\mu_1$.  But
$\mu_1\precneqq\nu_1$ implies that $S_{\mu_2}^* S_{\mu_1}^*
S_{\nu_1}=0$, a contradiction.
Suppose that $\mu_1$ and $\nu_1$ are of type $E_i\times F_i$.
If $\nu_1\npreceq\mu_1$ then
by \thmjj, we must have that $\mu_1$ and $\nu_1$ are separately
comparable in each coordinate, and that $\nu_1$ properly extends $\mu_1$
in at least one of the coordinates.  E.g. suppose
$S_{\mu_1}^*S_{\nu_1}=S_{(p,t(q)}S_{(t(p),q)}^*$ with $\ell(p)>0$. 
Then $S_{\mu_2}^*S_{\mu_1}^*S_{\nu_1} = S_{\mu_2}^*S_{(p,t(q))}S_{(t(p),q)}^*
=0$, again contradicting the nonzero hypothesis.

Moreover, if $|\nu|>1$ we must have $\mu_1=\nu_1$.  For if not,
$S_{\mu_1}^*S_{\nu_1}=S_{\mu_1'}^*$ where $\mu_1'$ is of the same type
as $\nu_1$, and $\ell(\mu_1')\not=(0,0)$.  But then
$S_{\mu_1'}^*S_{\nu_2}=0$, since $\nu_1$ and $\nu_2$ are of different
types.

Therefore if $|\nu|=1$ we have $\nu=\nu_1\preceq\mu_1\preceq\mu$,
while if $|\nu|>1$ we have $S_\mu^*S_\nu=S_{\mu_2\cdots\mu_n}^*
S_{\nu_2\cdots\nu_s}$, and the inductive hypothesis finishes the
argument.  \qed
\enddemo

We give two corollaries that will be needed in the proof of \thmww\
below.

\proclaim{\thmqq}  Let $\mu$, $\nu\in\Omega^*$ with
$S_\mu^*S_\nu\not=0$.  Suppose in addition that
$\ell(\mu)\ge\ell(\nu)$.  Then $\nu\preceq\mu$.
\endproclaim

\proclaim{\thmrr}  Let $\mu$, $\nu\in\Omega^*$ with
$S_\mu^*S_\nu\not=0$.  Suppose that $\ell(\mu)=(j,k)$ with $j<k$, and 
$\ell(\nu)\le(k,k)$.  Then either $\nu\preceq\mu$, or there are $p\in 
E_i^*$ and $q\in F_i^*$, with $\ell(p)>0$, such that $S_\mu^*S_\nu = 
S_{(p,t(q)}S_{(t(p),q)}^*$.  An analogous result holds with the roles 
of the two coordinates in $E_i\times F_i$ reversed.
\endproclaim

It follows from \thmpp\ that 
$\Theta$ is spanned by elements of the form $S_\mu S_\nu^*$ for which 
$t(\mu)=t(\nu)$. 
It also follows from the relations that there is an action, 
$\beta$, of $\t^2$ on $\Theta$, defined by
$\beta_z(S_\mu S_\nu^*)=\zeta_{\ell(\mu)-\ell(\nu)}(z)\,S_\mu S_\nu^*$.
We note that $\pi:\Theta\to A$ is equivariant for $\beta$ and $\alpha$.
We make some elementary computations in $\Theta\times_\beta \t^2$.
By means of the surjection $\pi$ we see that the analogous results 
hold in $A\times_\alpha\t^2$.
The elements (in $C(\t^2,\Theta)\subseteq \Theta\times_\beta\t^2$)
of the form $\zeta_mS_\mu S_\nu^*$, where $\mu$, $\nu\in\Omega^*$, 
make up a total subset of $\Theta\times_\beta\t^2$.  The fixed-point 
algebra $\Theta^\beta$ sits inside $\Theta\times_\beta\t^2$ as the closed 
linear span of the constant functions $S_\mu S_\nu^*$ for which 
$\ell(\mu)=\ell(\nu)$.  We recall the formulas for multiplication and 
adjoint in $C(\t^2,\Theta)\subseteq \Theta\times_\beta\t^2$:
$$\align
f\cdot g(z)&=\int f(v)\beta_v\bigl(g(v\inv z)\bigr)\,dv, \\
f^*(z)&=\beta_z\bigl(f(z\inv)^*\bigr). \\
\endalign$$

\definition{\thmnn} For $i=0$, 1, 
let $E_{i,1}\subseteq E_{i,2}\subseteq \cdots$ 
and $F_{i,1}\subseteq F_{i,2}\subseteq \cdots$ be finite irreducible 
non-circuit subgraphs of $E_i$ and $F_i$, with $\bigcup_k E_{i,k} 
=E_i$ and $\bigcup_k F_{i,k}=F_i$.
Let $\Omega_k$ denote the 
subobject of $\Omega$ comprising $D$,  $E_{0,k}\times 
F_{0,k}$, and $E_{1,k}\times F_{1,k}$.  We further let
$$\align
X_k&=\bigcup_{i,j\le k}\Omega_k^{(i,j)} \\
\intertext{and}
\Theta^\beta_k&=\text{span}\,\{S_\mu S_\nu^* : \mu,\nu\in X_k,\
\ell(\mu)=\ell(\nu)\}.
\endalign$$
\enddefinition

\proclaim{\thmzz} $\Theta^\beta$ is an AF algebra, with
$\{\Theta^\beta_k : k=1,2,\ldots\}$ as an approximating system
of finite dimensional \cstar-subalgebras.
\endproclaim

\demo{Proof}
It is clear that $\{S_\mu S_\nu^* : \mu,\nu\in X_k,\
\ell(\mu)=\ell(\nu)\}$ is a finite set, and
it follows from \thmpp\ that $\Theta^\beta_k$ is a finite dimensional 
\cstar-algebra.  Since $\cup_k\Theta^\beta_k$ is dense in
$\Theta^\beta$, $\Theta^\beta$ is an AF algebra.
\enddemo

\proclaim{\thmxx} $A^\alpha =  C^*\bigl(c\inv(0)\bigr)$ is an AF 
algebra.
\endproclaim

\demo{Proof}  Since $\pi:\Theta\to A$ is equivariant, and 
$\Theta^\alpha$ is AF, then so is $A^\alpha$.
For any $r$-discrete groupoid $G$ with continuous cocycle $c$ taking 
values in a discrete abelian group, and $\alpha$ the induced action of 
the dual group on $C^*(G)$,
it is a fact that 
$C^*(G)^\alpha = C^*\bigl(c\inv(0)\bigr)$.  \qed
\enddemo

A short computation shows that $\Theta^\beta$ is a 
hereditary subalgebra of $\Theta\times_\beta\t^2$.
The definition of ${\widehat\beta}$ is:
${\widehat\beta}_n(f)(z)=\zeta_n(z)\,f(z)$, for $f\in C(\t^2,\Theta)$, 
from which we find
$${\widehat\beta}_n(\zeta_m S_\mu S_\nu^*)=\zeta_{m+n}S_\mu S_\nu^*.$$

\proclaim\thmgg Let $\sigma$, $\sigma'$, $\mu$, and $\nu$ be finite paths with
$t(\mu)=t(\nu)$, $t(\sigma)=o(\mu)$, and $t(\sigma')=o(\nu)$, and such that
$\ell(\sigma\mu)=\ell(\sigma'\nu)$.  
Then $\zeta_{\ell(\sigma')}S_\mu S_\nu^* = S_\sigma^*\cdot 
S_{\sigma\mu}S_{\sigma'\nu}^* \cdot \zeta_{\ell(\sigma')}S_{\sigma'}$,
where $\cdot$ represents 
multiplication in $\Theta\times_\beta\t^2$.
\endproclaim

\demo{Proof} Let $m=\ell(\sigma')$.  We compute:
$$\align
S_\sigma^*\cdot S_{\sigma\mu}S_{\sigma'\nu}^*\cdot 
\zeta_mS_{\sigma'}(z)&=
\int\int S_\sigma^*\beta_{v}\Bigl(S_{\sigma\mu}S_{\sigma'\nu}^* 
\beta_{w}\bigl( 
(\zeta_mS_{\sigma'})(w\inv v\inv z)\bigr)\Bigr)\,dw\,dv \\
&=\int\int S_\sigma^*S_{\sigma\mu}S_{\sigma'\nu}^*\beta_{vw}
(\zeta_m(w\inv v\inv z)S_{\sigma'})\, dw\,dv \\
&=S_\mu S_{\sigma'\nu}^*\int\int\zeta_m(w\inv v\inv 
z)\zeta_m(vw)S_{\sigma'}\, 
dw\,dv \\
&=\zeta_{\ell(\sigma')}(z)S_\mu S_\nu^*. \qed
\endalign$$
\enddemo

Let $I_0$ denote the ideal in $\Theta\times_\beta\t^2$ generated by 
$\Theta^\beta$.  Note that it follows from \thmgg\ that 
$\zeta_{\ell(\sigma')}S_\mu S_\nu^*\in I_0$, 
where $\mu$, $\nu$, and $\sigma'$ are as in the statement.

\proclaim{\thmz} ${\widehat\beta}_{(1,1)}(I_0)\subseteq I_0$.
\endproclaim

\demo{Proof} It suffices to show that ${\widehat\beta}_{(1,1)}(\Theta^\beta) 
\subseteq I_0$.  So let $\mu$ and $\nu$ be finite paths with 
$t(\mu)=t(\nu)$ and $\ell(\mu)=\ell(\nu)$.  
Let $\sigma$ and $\sigma'$ be  finite paths with 
$\ell(\sigma)=\ell(\sigma')=(1,1)$, $t(\sigma)=o(\mu)$, and 
$t(\sigma')=o(\nu)$ 
($\sigma$ and $\sigma'$ exists since edges in $D$ 
map to $(1,1)$ under $c$).  Then we have
$$\align
{\widehat\beta}_{(1,1)}(S_\mu S_\nu^*)&= 
\zeta_{(1,1)}S_\mu S_\nu^* \\
&=S_\sigma^*\cdot S_{\sigma\mu}S_{\sigma'\nu}^*\cdot \zeta_{(1,1)}S_{\sigma'},
\text{ by \thmgg,} \\
&\in I_0. \qed
\endalign$$
\enddemo

We will let $I_n=({\widehat\beta}_{(1,1)})^n(I_0)$ for $n\in\z$.  We have a 
composition series:
$$\cdots\vartriangleleft I_1\vartriangleleft I_0\vartriangleleft 
I_{-1}\vartriangleleft \cdots \vartriangleleft \Theta\times_\beta\t^2.$$

\proclaim{\thmhh} 
$\Theta\times_\beta\t^2 = \overline{\bigcup_{n\in\z}I_n}$.
\endproclaim

\demo{Proof} Let $m\in\z^2$ and finite paths $\mu$, $\nu$ with 
$t(\mu)=t(\nu)$ be given.  We will show that there exists $k\in\z$ such 
that $\zeta_mS_\mu S_\nu^*\in I_{-k}$.  Choose $k\in\z$ with
$$\align
m+(k,k)&\ge(1,1),\text{ and }\tag\displayj \\
m+(k,k)+\ell(\nu)&\ge \ell(\mu)+(1,1). \tag\displayk
\endalign$$
It follows from (\displayj)
that there exists a finite path $\sigma'$ with 
$t(\sigma')=o(\nu)$ and $\ell(\sigma')=m+(k,k)$ 
(this is because every vertex is the 
terminus of a path of length $(1,1)$ with origin in $E_i\times F_i$ for 
some $i$).  Then it follows from (\displayk)
that there exists a finite 
path $\sigma$ with $t(\sigma)=o(\mu)$ and $\ell(\sigma)=\ell(\sigma')
-\ell(\mu)+\ell(\nu)$.  We have
$$\align
({\widehat\beta}_{(1,1)})^k(\zeta_mS_\mu S_\nu^*)&= 
\zeta_{m+(k,k)}S_\mu S_\nu^* \\
&=S_\sigma^*\cdot S_{\sigma\mu}S_{\sigma'\nu}^*\cdot 
\zeta_{m+(1,1)}S_{\sigma'},
\text{ by \thmgg,} \\
&\in I_0.
\endalign$$
It follows that $\zeta_m S_\mu S_\nu^*\in I_{-k}$. \qed
\enddemo

We now wish to
show that $\pi:\Theta\to A$ is an isomorphism.  To accomplish
this we need a detailed study of the finite dimensional approximating 
subalgebras of $\Theta^\beta$.  (This analysis will be needed again in
\thmbb.)  We begin with several definitions.  We remark that in
general, the 
structure of the finite dimensional subalgebras of the AF core of a
Toeplitz graph algebra
is made complicated 
by the fact that the Cuntz-Krieger relations are not satisfied at all 
vertices, and hence that there will be nonzero {\it defect\/}
projections at some vertices.  
In \thmyy\ below we treat this situation, which is
analogous to, but much easier, than the situation of $\Theta^\beta$.
The case of $\Theta^\beta$ is further complicated by the hybrid graph 
structure of $\Omega$.

\definition{\thmss} For $x\in\Omega_k^{(0,0)}$ we define projections
$\lambda_k(x)$, $\rho_k(x)$, $\omega_k^{(0)}(x)$ and
$\omega_k(x)\in\Theta^\beta_k$ by
$$\aligned
\lambda_k(x)&=P_x - \sum\{S_\mu S_\mu^* : \mu\in\Omega_k^{(1,0)}(x)\} 
\\
\rho_k(x)&=P_x - \sum\{S_\mu S_\mu^* : \mu\in\Omega_k^{(0,1)}(x)\} 
\\
\omega_k(x)&=\lambda_k(x)\rho_k(x) \\
\omega_k^{(0)}(x)&=\cases \omega_k(x),&\text{if
$x\not\in\{u_0,u_1\}$} \\
\omega_k(x)-S_{\alpha_i}S_{\alpha_i}^* -
S_{\epsilon_i}S_{\epsilon_i}^*,&\text{if $x=u_i$}.\endcases
\endaligned$$
\enddefinition

\remark{\thmtt}
Note that $\omega_k^{(0)}(u_i)$ is a projection by 
relation (iv) of \thmjj. 
We list some easy consequences of the definitions.  Recall (from the
end of section 1) that $S(E_{i,k})$ is the set of vertices of
$E_{i,k}$ that do not emit edges in $E_i$ other than those emitted in 
$E_{i,k}$.
$$\align
\lambda_k(x)&\not=0 \text{ iff } x=(y,z) \text{ with } y\not\in
S(E_{i,k}). \\
\rho_k(x)&\not=0 \text{ iff } x=(y,z) \text{ with } z\not\in
S(F_{i,k}). \\
\omega_k(x)&\not=0 \text{ iff } x=(y,z) \text{ with } y\not\in
S(E_{i,k}) \text{ and } z\not\in
S(F_{i,k}). \\
\lambda_k(x)S_\mu&=0 \text{ whenever } \ell(\mu)\ge(1,0),\
\mu\in\Omega_k^*. \\
\rho_k(x)S_\mu&=0 \text{ whenever } \ell(\mu)\ge(0,1),\
\mu\in\Omega_k^*. \\
\omega_k(x)S_\mu&=0 \text{ whenever } \ell(\mu)\not=(0,0),\
\mu\in\Omega_k^*. \\
\lambda_k\bigl(o(\mu)\bigr)S_\mu &=S_\mu\lambda_k \bigl(t(\mu)\bigr)
\text{ if }\ell(\mu)=(0,1),\ \mu\in\Omega_k^*. \\
\rho_k\bigl(o(\mu)\bigr)S_\mu &=S_\mu\rho_k \bigl(t(\mu)\bigr)
\text{ if }\ell(\mu)=(1,0),\ \mu\in\Omega_k^*.
\endalign$$
\endremark

\definition{\thmuu}  $\mu\in X_k$ is {\it maximal} if whenever
$\mu\preceq\mu'$ with $\mu'\in X_k$ then $\mu'=\mu$.
\enddefinition

We remark that $\mu\in X_k$ is maximal if and only if one of the following
occurs:
\roster
\item $\ell(\mu)=(k,k)$.
\item $\ell(\mu)=(k-1,k)$ or $(k,k-1)$, and $t(\mu)\in\{a_0,a_1\}$.
\endroster

\definition{\thmvv}  We define certain non-zero projections in 
$\Theta^\beta_k$.  They are of four kinds:
\roster
\item $S_\mu S_\mu^*$, where $\mu\in X_k$ is maximal.
\smallskip
\item $S_\mu\lambda_k\bigl(t(\mu)\bigr)S_\mu^*$, where $\ell(\mu)=(j,k)$ 
with $j<k$, and $t(\mu)=(y,z)$ with $y\not\in S(E_{i,k})$.
\smallskip
\item $S_\mu\rho_k\bigl(t(\mu)\bigr)S_\mu^*$, where $\ell(\mu)=(k,j)$ 
with $j<k$, and $t(\mu)=(y,z)$ with $z\not\in S(F_{i,k})$.
\smallskip
\item $S_\mu\omega_k\bigl(t(\mu)\bigr)S_\mu^*$, where 
$\ell(\mu)\le(k-1,k-1)$ 
with $j<k$, and $t(\mu)=(y,z)$ with $y\not\in S(E_{i,k})$ and
$z\not\in S(F_{i,k})$.
\endroster
We let $\theta_i(\mu)$ denote the projection in \thmvv($i$) above
constructed from the path $\mu$.
\enddefinition

\proclaim{\thmww}  
\roster
\runinitem $\theta_i(\mu)$ is a minimal projection in $\Theta^\beta_k$.
\smallskip
\item $\theta_i(\mu)$ and $\theta_{i'}(\nu)$ are equivalent in
$\Theta^\beta_k$ if and only if $i=i'$, $\ell(\mu)=\ell(\nu)$, and
$t(\mu)=t(\nu)$.
\smallskip
\item $\sum_{i,\mu}\theta_i(\mu)=1_{\Theta^\beta_k}$.
\endroster
\endproclaim

\demo{Proof} First note that if $\ell(\mu)=\ell(\nu)$ and
$t(\mu)=t(\nu)$ then $(S_\nu S_\mu^*)\theta_i(\mu)(S_\nu S_\mu^*)^* = 
\theta_i(\nu)$, proving the reverse direction of (2).  We  now claim
that if $\sigma$, $\tau\in X_k$ with $\ell(\sigma)=\ell(\tau)$
are such that $\theta_i(\mu)S_\sigma S_\tau^*\theta_{i'}(\nu)\not=0$,
then $i=i'$, $t(\mu)=t(\nu)$, $\ell(\mu)=\ell(\nu)$, and moreover,
$\sigma\preceq\mu$, $\tau\preceq\nu$.  This will conclude the proof of
(2).  (1) will also follow, since then $\theta_i(\mu)S_\sigma
S_\tau^*\theta_i(\mu)\not=0$ implies that $\sigma$, $\tau\preceq\mu$. 
Then $\ell(\sigma)=\ell(\tau)$ implies that $\sigma=\tau$, and hence
that $\theta_i(\mu)S_\sigma S_\tau^*\theta_i(\mu)=\theta_i(\mu)$.

To prove the claim, assume that
$\theta_i(\mu)S_\sigma S_\tau^*\theta_{i'}(\nu)\not=0$.  
We first consider the situation
$S_\mu^*S_\sigma\not=0$.  Suppose first that $\mu$ is maximal.  If
$\ell(\mu)=(k,k)$ then $\sigma\preceq\mu$ by \thmqq.  If
$\ell(\mu)\not=(k,k)$ then $t(\mu)\in\{a_0,a_1\}$.  By \thmpp,
$t(\sigma)=t(\mu)$.  Since $\mu$ is maximal, $\mu\precneqq\sigma$ is
impossible.  Then by \thmpp\ we have $\sigma\preceq\mu$.  Next suppose
that $\ell(\mu)=(j,k)$ with $j<k$ and $t(\mu)=(y,z)$ with $y\not\in
S(E_{i,k})$.  If $\sigma\npreceq\mu$, \thmrr\ implies that
$S_\mu^*S_\sigma=S_{(p,t(q))}S_{(t(p),q)}^*$ with $\ell(p)>0$.  But
then
$\lambda\bigl(t(\mu)\bigr)S_\mu^*S_\sigma=\lambda\bigl(t(\mu)\bigr)
S_{(p,t(q))}S_{(t(p),q)}^*=0$.  Hence we must have $\sigma\preceq\mu$.
An analogous argument handles the case 
$\ell(\mu)=(k,j)$ with $j<k$ and $t(\mu)=(y,z)$ with $z\not\in
S(F_{i,k})$.  Finally suppose $\ell(\mu)\le(k-1,k-1)$ and
$t(\mu)=(y,z)$ with $y\not\in S(E_{i,k})$ and $z\not\in S(F_{i,k})$.  
If $\mu\precneqq\sigma$ then $\sigma=\mu\sigma'$ with
$\ell(\sigma')\not=(0,0)$.  Then
$$\omega_k(y,z)S_\mu^*S_\sigma=\omega_k(y,z)S_{\sigma'}=0,$$
a contradiction.  If $\sigma\npreceq\mu$, then \thmpp\ implies that
$S_\mu^*S_\sigma=S_{(p,t(q))}S_{(t(p),q)}^*$ with $\ell(p)\not=0$ and 
$\ell(q)\not=0$.  But then
$$\omega_k(y,z)S_\mu^*S_\sigma=\omega_k(y,z)
S_{(p,t(q))}S_{(t(p),q)}^*=0,$$
a contradiction.  Thus in all cases we have $\sigma\preceq\mu$.

Now write $\mu=\sigma\mu'$.  We have $S_\mu^*S_\sigma S_\tau^* =
S_{\mu'}^*S_\tau^*=S_{\tau\mu'}^*$.  Note that $t(\mu)=t(\tau\mu')$
and $\ell(\mu)=\ell(\tau\mu')$.  We then have
$$\align
\theta_i(\mu)S_\sigma S_\tau^*&= \theta_i(\mu)S_\mu S_\mu^* S_\sigma
S_\tau^* \\
&=\theta_i(\mu)S_\mu S_{\tau\mu'}^* \\
&=S_\mu S_{\tau\mu'}^*\theta_i(\tau\mu') \\
\theta_i(\mu)S_\sigma S_\tau^*\theta_{i'}(\nu)&= S_\mu
S_{\tau\mu'}^*\theta_i(\tau\mu')\theta_{i'}(\nu).
\endalign$$
Thus $\theta_i(\tau\mu')\theta_{i'}(\nu)\not=0$, and hence
$\theta_i(\tau\mu')S_\sigma S_\sigma^*\theta_{i'}(\nu)\not=0$
for $\sigma=\tau\mu'$ and for $\sigma=\nu$.  It follows from the above that
$\tau\mu'\preceq\nu$ and $\nu\preceq\tau\mu'$, 
proving the claim and finishing the proof of (1) and (2).

Before proving (3) we make some preliminary observations.  
Let $x\in\Omega_k^{(0,0)}$. Using
\thmss\ and \thmtt\ we find that
$$\align
P_x&=\lambda_k(x)+ \sum\{S_\mu S_\mu^* : \mu\in\Omega_k^{(1,0)}(x)\}
\tag\displaya \\
P_x&=\rho_k(x) + \sum\{S_\mu S_\mu^* : \mu\in\Omega_k^{(0,1)}(x)\}
\tag\displayb \\
P_x&=\left(\lambda_k(x)+ \sum\{S_\mu S_\mu^* :
\mu\in\Omega_k^{(1,0)}(x)\}\right) \\
&\qquad\left(
\rho_k(x) + \sum\{S_\mu S_\mu^* : \mu\in\Omega_k^{(0,1)}(x)\} \right) 
\\
&=\lambda_k(x)\rho_k(x) + 
\sum\{\lambda_k(x)S_\mu S_\mu^* : \mu\in\Omega_k^{(1,0)}(x)\} \\
&\qquad +\ \sum\{\rho_k(x)S_\nu S_\nu^* : \nu\in\Omega_k^{(0,1)}(x)\} \\
&\qquad\qquad +\ \sum\{S_{\mu\nu} S_{\mu\nu}^* : \mu\in\Omega_k^{(1,0)}(x),\
\nu\in\Omega_k^{(0,1)}(x)\} \\
&=\omega_k(x) + 
\sum\{S_\mu\lambda_k\bigl(t(\mu)\bigr) S_\mu^*
     : \mu\in\Omega_k^{(1,0)}(x)\}\tag\displayc \\
&\qquad +\ \sum\{S_\nu\rho_k\bigl(t(\nu)\bigr) S_\nu^*
     : \nu\in\Omega_k^{(0,1)}(x)\} \\
&\qquad\qquad +\ \sum\{S_{\mu} S_{\mu}^* : \mu\in\Omega_k^{(1,1)}(x)\}.
\endalign$$
Next, let $\mu_1\in\Omega_k^{(0,1)}$ and $\mu_2\in\Omega_k^{(1,0)}$.
If $t(\mu_1)=o(\mu_2)$ then
$$\lambda_k\bigl(o(\mu_1)\bigr) S_{\mu_1} \rho_k\bigl(o(\mu_2)\bigr)
S_{\mu_2} = S_{\mu_1} \lambda_k\bigl(t(\mu_1)\bigr) S_{\mu_2}
\rho_k\bigl(t(\mu_2)\bigr) =0,$$
and if $t(\mu_2)=o(\mu_1)$ then
$$\rho_k\bigl(o(\mu_2)\bigr) S_{\mu_2} \lambda_k\bigl(o(\mu_1)\bigr)
S_{\mu_1} = S_{\mu_2} \rho_k\bigl(t(\mu_2)\bigr) S_{\mu_1}
\lambda_k\bigl(t(\mu_1)\bigr) =0,$$
by \thmtt.  Finally, if $t(\mu_1)=u_i$ and
$\mu_3\in\{\alpha_i,\epsilon_i\}$ then
$$\lambda_k\bigl(o(\mu_1)\bigr) S_{\mu_1} S_{\mu_3} = S_{\mu_1}
\lambda_k(u_i) S_{\mu_3} = S_{\mu_1} S_{\mu_3},$$
by \thmss\ and \thmjj\ (iv).  Hence for
$\mu_3\in\Omega_k^{(1,1)}(u_i)$,
$$\lambda_k\bigl(o(\mu_1)\bigr) S_{\mu_1} S_{\mu_3} =
\cases S_{\mu_1\mu_3},\text{if }\mu_3\in\{\alpha_i,\epsilon_i\} \\
0,\text{otherwise}. \endcases$$
An analogous result holds with $\mu_2$ replacing $\mu_1$ and $\rho$
replacing $\lambda$.

Now we prove (3).  Let $x\in\Omega_k^{(0,0)}$, and consider the
expression in (\displayc).  In each term of the form
$S_\mu S_\mu^*= S_\mu P_{t(\mu)} S_\mu^*$, 
substitute for $P_{t(\mu)}$ with formula
(\displayc) if $\ell(\mu)\le(k-1,k-1)$, with formula (\displaya) if
$(0,k)\le\ell(\mu)\le(k-1,k)$, and with formula (\displayb) if
$(k,0)\le\ell(\mu)\le(k,k-1)$.  Use the above observations to
eliminate zero terms and to simplify, then repeat.  This process must 
stop in a finite number of steps, giving
$$\align
P_x&= \sum\{S_\mu \omega_k\bigl(t(\mu)\bigr) S_\mu^* : o(\mu)=x,\
\ell(\mu)\le(k-1,k-1)\} \\
&\qquad +\ \sum_{j=0}^{k-1}\sum\{ S_\mu \lambda_k\bigl(t(\mu)\bigr)
S_\mu^* : o(\mu)=x,\ \ell(\mu)=(j,k)\} \\
&\qquad\qquad +\ \sum_{j=0}^{k-1} \sum\{ S_\mu \rho_k\bigl(t(\mu)\bigr)
S_\mu^* : o(\mu)=x,\ \ell(\mu)=(k,j)\} \\
&\qquad\qquad\qquad +\ \sum\{S_\mu S_\mu^* : o(\mu)=x,\ \mu\in X_k \text{
maximal}\} \\
&=\sum_i \sum\{\theta_i(\mu) : o(\mu)=x\}.
\endalign$$
Now (3) follows by summing over $x\in\Omega_k^{(0,0)}$. \qed
\enddemo

\proclaim{\thmx} (Gauge-invariant uniqueness theorem)  Let
$$\phi:(\Theta,\t^2,\beta)\to(C,\t^2,\gamma)$$ 
be an equivariant
$*$-homomorphism between \cstar-dynamical systems.  If
$\phi\restrictedto{\Theta^\beta}$ is injective then $\phi$ is injective.
\endproclaim

\demo{Proof}  We have that 
$\widetilde\phi:\Theta\times_\beta\t^2\to C\times_\gamma\t^2$ is 
injective, since it is so on the ideals $I_n$.  Therefore 
$\widetilde{\widetilde\phi}:\Theta\times_\beta\t^2\times_{\widehat\beta}\z^2
\to C\times_\gamma\t^2\times_{\widehat\gamma}\z^2$ is injective.  The result 
now follows from Takesaki-Takai duality.  \qed
\enddemo

\proclaim{\thmy} $A=C^*(G)$ is isomorphic to $\Theta$, and
is simple, purely infinite,
nuclear and classifiable, i.e. a UCT Kirchberg algebra.
\endproclaim

\demo{Proof}  Let $\lambda:C^*(G)\to C^*_r(G)$ be the left regular
representation.  Then $\lambda\circ\pi:\Theta\to C^*_r(G)$.  The
action of $\t^2$ on $C^*(G)$ clearly descends equivariantly to
$C^*_r(G)$, so that $\lambda\circ\pi$ is equivariant. We note
that $\lambda\circ\pi\restrictedto{\Theta^\beta}$ is injective.  This 
follows from the facts that $\Theta^\beta$ is AF, and that the minimal
projections $\theta_i(\mu)$  (\thmww)
in the subalgebras $\Theta^\beta_k$ have
nonzero image in $C_c(G)^\alpha\subseteq C^*_r(G)^\alpha$. Thus
$\lambda\circ\pi$ is injective by \thmx.  It follows that $\lambda$ is
injective (this could also be deduced from nuclearity of $C^*(G)$,
proved below, and \bib\anantharaman).  We have that
$C^*(G)\times_\alpha\t^2=\widetilde\pi(\Theta\times_\beta \t^2)$ 
is AF, and hence $A$, which is strongly Morita equivalent to
$A\times_\alpha\t^2\times_\alhat\z^2$, is nuclear and classifiable.
Simplicity and pure infiniteness follow from \thmv. \qed
\enddemo

The following proposition is necessary for our application of the 
results of this paper in \bib\spielbergc.  The proof is immediate 
from the description of $\Theta$ by generators and relations.  (We 
remark that in that application, the graphs $F_i$ will be chosen to 
represent Kirchberg algebras having $K$-theory of the form $(\z,0)$ 
or $(0,\z)$, while the graphs $E_i$ will be chosen to represent 
Kirchberg algebras with preassigned $K_0$ and trivial $K_1$.)

\proclaim{\thmaaa}  Let $\Gamma_i$ be a subgroup of 
$\text{Aut}\,(E_i)$ fixing the vertex $v_i$.  There is a 
homomorphism $\Gamma_0\times\Gamma_1\to\text{Aut}\,(\Theta)$ defined 
on generators by letting $\Gamma_0$ act on $E_0$, $\Gamma_1$ on 
$E_1$, and the trivial action on $D$, $F_0$ and $F_1$.  Moreover, if 
$x_0\in E_0^0$ is fixed by $\Gamma_0$, then the corner of $\Theta$ 
defined by the projection $P_{(x_0,w_0)}$ is invariant for the action 
of $\Gamma_0\times \Gamma_1$.
\endproclaim

\head 4.  The $K$-theory of $A$ \endhead

We may now omit the use of the $*$-isomorphism $\pi$, and identify 
$\Theta$ with $A$.

It follows from the fact that the ideals $I_n$ are AF
that the inclusion of $I_n$ into 
$A\times_\alpha\t^2$ induces an injection in $K_0$.  We let $\phi$ 
denote the automorphism $\alhat_{(1,1)\,*}$ of 
$K_0(A\times_\alpha\t^2)$.  We have
$$K_0(A\times_\alpha\t^2)=\bigcup_{n\in\z}\phi^n\bigl(K_0(A^\alpha)\bigr).$$
We let $W=K_0(A\times_\alpha\t^2)$ and $W_n=K_0(I_n)$.  Thus 
$W_n\supseteq W_{n+1}$ and $W=\cup_{n\in\z}W_n$.  Since $W_0\cong 
K_0(A^\alpha)$, $W_0$ is generated by elements of the form 
$[S_\mu S_\mu^*]$, where $\mu$ is a finite path.  If $t(\mu)=t(\nu)$ and 
$\ell(\mu)=\ell(\nu)$ then $S_\mu S_\nu^*$
 is a partial isometry in $A^\alpha$ 
implementing an equivalence between $S_\mu S_\mu^*$ and $S_\nu S_\nu^*$, so 
that $[S_\mu S_\mu^*]=[S_\nu S_\nu^*]$.  

\proclaim{\thmaa} Let $\mu$ and $\sigma$ be finite paths with
$t(\sigma)=o(\mu)$. 
Then $\alhat_{\ell(\sigma)*}\bigl([S_\mu S_\mu^*]\bigr) = 
[S_{\sigma\mu}S_{\sigma\mu}^*]$.
\endproclaim

\demo{Proof} Let $f\in C(\t^2,A)$ be given by 
$f(z)=\zeta_{\ell(\sigma)}S_{\sigma\mu}S_\mu^*$.  Routine computations give
$$\align
f^*f&= \zeta_{\ell(\sigma)}S_\mu S_\mu^* \\
&=\alhat_{\ell(\sigma)}(S_\mu S_\mu^*),\text{ and} \\
ff^*&= S_{\sigma\mu}S_{\sigma\mu}^*. \qed
\endalign$$
\enddemo

We let $B$ denote the subalgebra of $A$ generated by the 
edges of $\Omega$ in $\cup_i E_i\times F_i$.  Thus $B$ is 
isomorphic to the direct sum of the algebras $C^*(E_i)\otimes 
C^*(F_i)$.  We note that $B$ is invariant under $\alpha$.  The 
crossed product $B\times_\alpha\t^2$ is an AF subalgebra of 
$A\times_\alpha\t^2$, and is isomorphic to a direct sum of
tensor products of AF 
algebras:  
$$B\times_\alpha\t^2\cong \bigoplus_i\; \bigl(C^*(E_i)\times\t\bigr) 
\otimes \bigl(C^*(F_i)\times\t\bigr).$$

\proclaim{\thmbb} Let $i:B\times_\alpha\t^2 \to A\times_\alpha\t^2$ be 
the inclusion map.  Then $i_*$ is injective in $K_0$.
\endproclaim

Before giving the proof, we give a preliminary lemma describing the 
finite dimensional approximants to the AF core of the
relative Toeplitz algebra of an ordinary graph.  We remark that this
is an easier version of \thmww.

\proclaim{\thmyy} Let $E$ be a finite directed graph.  
Let $S\subseteq E^0$ not contain any sink of $E$.  
For $k\ge 1$ 
let $C_k(E,S)$ be the finite dimensional \cstar-subalgebra of 
$\calt\oh(E,S)$ given by
$$\align
C_k(E,S)&=\text{\rm span}\,\{S_p S_q^* : p,\,q\in\cup_{j\le k}E^j,\
\ell(p)=\ell(q)\}. \\
\intertext{For $y\in E^0\setminus S$ let $\xi_y\in
C_k(E,S)$ be given by}
\xi_y&=P_y-\sum_{e\in E^1(y)} S_e S_e^*. \\
\intertext{For $0\le j<k$ and $y\in E^0\setminus S$ let}
N_j^{(k)}(y)&= \{S_p\xi_y S_p^* : p\in E^j,\ t(p)=y\}. \\
\intertext{For $y\in E^0$\ let}
N_k^{(k)}(y)&= \{S_p S_p^* : p\in E^k,\ t(p)=y\}. \\
\intertext{Set}
N^{(k)}&= \bigcup\bigl\{ N_j^{(k)}(y) : 0\le j<k,\ y\in E^0\setminus 
S\bigr\} \cup \bigcup\bigl\{ N_k^{(k)}(y) : y\in E^0\bigr\}.
\endalign$$
Then
\roster
\item"(\displayd)" $N^{(k)}$ is a basis for 
$C_k(E,S)$ consisting of pairwise 
orthogonal minimal projections.
\item"(\displaye)" Let $a\in N_i^{(k)}(y)$ and 
$b\in N_j^{(k)}(z)$, where $0\le 
i,j\le k$ and $y$, $z\in 
E^0$.  Then $a$ and $b$ are equivalent in 
$C_k(E,S)$ if and only if $i=j$ and $y=z$.
\endroster
\endproclaim
\demo{Proof}  We first note that
$$\align
\xi_y S_e &=0 \text{ for } y\in E^0\setminus S,\ e\in E^1, 
\tag\displayf \\
S_p^* S_q&=0 \text{ iff } p \text{ and } q \text{ are not comparable,
for } p,\ q\in E^*. \tag\displayg
\endalign$$
It follows easily  that
$$\align
S_p \xi_{t(p)} S_p^* S_q \xi_{t(q)} S_q^* &\not=0
\text{ iff } p=q \text{ and } t(p),\
t(q)\in E^0\setminus S; \\
S_p S_p^* S_q S_q^* &=0 \text{ if } \ell(p)=\ell(q)=k \text{ and }
p\not=q; \\
S_p \xi_{t(p)} S_p^* S_q S_q^* &=0 \text{ if } \ell(p)<\ell(q) \text{ 
and } t(p)\in E^0\setminus S.
\endalign$$
Therefore the projections of $N^{(k)}$ are pairwise orthogonal.  We now
consider the minimality and equivalence of projections in $N^{(k)}$
together.  For this, fix paths $r$, $s\in\cup_{j\le k}E^j$
with $\ell(r)=\ell(s)$.  Let 
$t(p)$, $t(q)\in E^0\setminus S$ with $\ell(p)$, $\ell(q)<k$, and
suppose that
$$S_p \xi_{t(p)} S_p^* S_r S_s^* S_q \xi_{t(q)} S_q^* \not=0.
\tag\displayh$$
Then (\displayf) and (\displayg) imply that $r\preceq p$, $s\preceq q$, and
$S_p^* S_r = S_q^* S_s$.  Hence $t(p)=t(q)$ and $\ell(p)=\ell(q)$, 
as required by (\displaye).  Moreover if $p=q$ then the product in
(\displayh) equals $S_p \xi_{t(p)} S_p^*$, proving that $S_p \xi_{t(p)}
S_p^*$ is minimal.  
If $p$, $q\in E^k$ and if
$$S_p S_p^* S_r S_s^* S_q S_q^*\not=0, \tag\displayi$$
then $r\preceq p$, $s\preceq q$, and $S_r^* S_p = S_s^* S_q$, 
so $t(p)=t(q)$ as required by (\displaye).  Moreover if $p=q$ then the
product in (\displayi) equals $S_p S_p^*$, proving that $S_p S_p^*$ is
minimal.  Finally, let $t(p)\in E^0\setminus S$ with $\ell(p)<k$,
and $\ell(q)=k$, and consider
$S_p\xi_{t(p)}S_p^*S_rS_s^*S_qS_q^*$.  It follows from 
(\displayf) and (\displayg) that $r\preceq p$ and $s\preceq q$.  
Then $S_rS_s^*S_q=S_{q'}$ where $\ell(q')=k$.  Since $\ell(p)<k$, it 
follows from the same considerations that $\xi_{t(p)}S_p^*S_{q'}=0$.
Thus $S_p\xi_{t(p)}S_p^*$ and $S_qS_q^*$ are inequivalent. 

For the reverse implication in (\displaye), note that if
$\ell(p)=\ell(q)<k$ and $y=t(p)=t(q)\in E^0\setminus S$ then
$$\align
(S_q \xi_y S_p^*)^* (S_q \xi_y S_p^*)&=S_p \xi_y S_p^* \\
(S_q \xi_y S_p^*) (S_q \xi_y S_p^*)^*&=S_q \xi_y S_q^*, \\
\intertext{and if $\ell(p)=\ell(q)=k$ and $t(p)=t(q),$ then}
(S_q S_p^*)^*(S_q S_p^*)&=S_p S_p^* \\
(S_q S_p^*)(S_q S_p^*)^*&=S_q S_q^*.
\endalign$$

Finally, we show that $\sum N^{(k)}=1$.  For convenience we will let
$\xi_y=0$ for $y\in S$.  Then for any $y\in E^0$ we have
$$\align
P_y&= \xi_y + \sum_{e_1\in E^1(y)} S_{e_1} S_{e_1}^* \\
&= \xi_y + \sum_{e_1\in E^1(y)} S_{e_1}\Bigl( \xi_{t(e_1)} +
\sum_{e_2\in E^1\bigl(t(e_1)\bigr)} S_{e_2} S_{e_2}^*\Bigr)
S_{e_1}^* \\
&=\xi_y + \sum_{p\in E^1(y)} S_p \xi_{(t(p)} S_p^* + \sum_{p\in E^2(y)}
S_p S_p^* \\
&=\cdots \\
&= \sum_{i=0}^{k-1} \sum_{p\in E^i(y)} S_p \xi_{t(p)} S_p^* + 
\sum_{p\in E^k(y)} S_p S_p^*.
\endalign$$
The result follows by summing over $y\in E^0$. \qed
\enddemo

\demo{Proof of \thmbb}  By repeated application of $\phi$ it suffices to show 
that $i_*:K_0(B^\alpha)\to K_0(A^\alpha)$ is injective.  Letting
$A^\alpha_k=\Theta^\beta_k$ (via the isomorphism $\pi$), we have
$A^\alpha=\overline{\cup_k A^\alpha_k}$.  We define finite dimensional
approximating subalgebras to $B^\alpha$ in a manner similar to the
$\Theta^\beta_k$.  Namely, let
$$X_k^{(0)}= \cup_{i=0}^1 \cup_{j,j'\le k} E_{i,k}^j \times
F_{i,k}^{j'}.$$
(Thus $X_k^{(0)}$ is the set of paths $\mu$ in $X_k$ that do not contain
edges from $D$.)  Now let
$$B^\alpha_k=\text{span}\,\bigl\{ S_\mu S_\nu^* : \mu,\ \nu\in
X_k^{(0)},\ \ell(\mu)=\ell(\nu)\bigr\}.$$
It is clear that $B^\alpha=\overline{\cup_k B^\alpha_k}$ and that
$B^\alpha_k\subseteq A^\alpha_k$.  We will show that the inclusion 
$B^\alpha_k\subseteq A^\alpha_k$ induces an injection in $K_0$.  This 
will prove the lemma.

We note that 
$$B^\alpha_k\cong\bigoplus_i C_k\bigl(E_{i,k},S(E_{i,k})\bigr) \otimes
C_k\bigl(F_{i,k},S(F_{i,k})\bigr).$$
Thus every minimal projection in $B^\alpha_k$ is a tensor product of
minimal projections from $C_k\bigl(E_{i,k},S(E_{i,k})\bigr)$ and
$C_k\bigl(F_{i,k},S(F_{i,k})\bigr)$, and two such are equivalent in
$B^\alpha_k$ if and only if they are separately equivalent in each
factor.  Note that if $(y,z)\in \Omega_k^{(0,0)}$ then
$$\align
\lambda_k(x)&= \xi_y\otimes P_z \\
\rho_k(x)&= P_y\otimes \xi_z \\
\omega_k^{(0)}(x)&= \xi_y\otimes\xi_z, \\
\intertext{while if $\mu=(p,q)\in X_k^{(0)}$, then}
S_\mu S_\mu^*&= S_p S_p^*\otimes S_q S_q^*.
\endalign$$
For $\mu\in X_k^{(0)}$ with $\ell(\mu)\le(k-1,k-1)$ and $t(\mu)=(y,z)$ 
with $y\not\in S(E_{i,k})$, $z\not\in S(F_{i,k})$, we define
$$\theta_4^{(0)}(\mu)=S_\mu \omega_k^{(0)}\bigl(t(\mu)\bigr) S_\mu^* 
\in B^\alpha_k,$$
analogously to \thmvv(4).  Then the projections in $B^\alpha_k$ of the
form $\theta_1(\mu)$, $\theta_2(\mu)$, $\theta_3(\mu)$ and
$\theta_4^{(0)}(\mu)$ form a complete family of pairwise orthogonal
minimal projections.  By \thmyy\ we see that they satisfy the
conditions for equivalence given in \thmww(2).  From \thmss\ we see
that
$$\theta_4^{(0)}(\mu)=\theta_4(\mu) + \theta_{i_1}(\tau_1) +
\theta_{i_2}(\tau_2) + \cdots,$$
where $\ell(\tau_1)$, $\ell(\tau_2)$, $\ldots > \ell(\mu)$.
This observation has the following 
consequence.  Choose bases for $K_0(A^\alpha_k)$ and $K_0(B^\alpha_k)$ 
consisting of classes of minimal projections as above.  If the bases 
are ordered by increasing length of the underlying paths, then the 
matrix of the map $K_0(B^\alpha_k)\to K_0(A^\alpha_k)$ induced from 
inclusion is lower triangular, with 1's on the diagonal.  Thus the map 
is injective.  \qed
\enddemo

We let $Y$ denote $K_0(B\times_\alpha\t^2)$.
By \thmbb\ we may identify $Y$ with $i_*(Y)\subseteq W$.  We now give 
a key lemma, that is based on the fact that the (ordinary) graph 
$D$ connecting the 2-graphs $E_i\times F_i$ is a ``bit of $\oh_2$''.

\proclaim{\thmcc} $(2\phi-\id)W\subseteq Y$.
\endproclaim

\demo{Proof} Let $x\in W$.  For $n\in\z$ large enough we have 
$\phi^n(x)\in W_0$.  If $(2\phi-\id)\phi^n(x)\in Y$, then since $Y$ 
is $\phi$-invariant we get $(2\phi-\id)(x)\in\phi^{-n}Y=Y$.  So we may 
assume that $x\in W_0$.  Since $W_0=K_0(A^\alpha)$ is generated by 
elements of the form $[S_\mu S_\mu^*]$ for finite paths $\mu$, we may assume 
that $x=[S_\mu S_\mu^*]$.

\demo{Case (i)} Suppose $t(\mu)\in E^0_i\times  F^0_i$ for some $i$.  
Since the $K_0$-class of the projection is unchanged if the path is 
replaced by a new path with the same length and terminus, we may 
assume that $\mu\in (E_i\times F_i)^*$, and so that $x\in Y$.
\enddemo

\demo{Case (ii)}  Suppose $t(\mu)\in \{a_0,a_1\}$.  
For definiteness we suppose $t(\mu)=a_0$.  For the rest of this
argument we will omit the subscript on $a_0$, $\beta_0$, $\gamma_0$ and
$\delta_0$.  Then 
we may assume that $\mu=\nu\beta^m$ for some $m\ge0$ and some path $\nu$ 
with $t(\nu)=a$.  We note that
$$\align
S_{\beta^m}S_{\beta^m}^*&=S_{\beta^m}(S_\beta S_\beta^*+S_\gamma 
S_\gamma^* + S_\delta S_\delta^*)S_{\beta^m}^*, \\
\intertext{while}
\phi[S_{\beta^m}S_{\beta^m}^*]&=[S_{\beta^{m+1}}S_{\beta^{m+1}}^*]. \\
\intertext{Hence}
[S_{\beta^m}S_{\beta^m}^*]&=2[S_{\beta^{m+1}}S_{\beta^{m+1}}^*] + 
[S_{\beta^m\delta}S_{\beta^m\delta}^*] \\
&\in2\phi[S_{\beta^m}S_{\beta^m}^*]+Y,
\endalign$$
since $t(\beta^m\delta)=u_i$.  Therefore 
$(2\phi-\id)[S_{\beta^m}S_{\beta^m}^*]\in Y$.  Thus
$$\align
(2\phi-\id)[S_\mu S_\mu^*]&=
(2\phi-\id)\phi[S_{\nu\beta^m}S_{\nu\beta^m}^*] \\
&=\alhat_{\ell(\nu)*}\circ(2\phi-\id)[S_{\beta^m}S_{\beta^m}^*],\text{ by 
\thmaa,} \\
&\in\alhat_{\ell(\nu)*}(Y) \\
&\subseteq Y. \qed
\endalign$$
\enddemo
\enddemo

\proclaim{\thmdd} $\ker\,(\id-\phi)\subseteq Y$.
\endproclaim

\demo{Proof} Let $x\in\ker\,(\id-\phi)$.  Then
$$\align
x&=\phi(x) \\
&=\phi(x)-(\id-\phi)(x) \\
&=(2\phi-\id)(x) \\
&\in Y. \qed
\endalign$$
\enddemo

The preceding and following lemmas will allow us to show that
the $K$-theory of $A\times_\alpha\t^2\times_\phi\z$ is given by the 
subalgebra $B$.  We let $\psi=\alhat_{(1,0)*}$, so that $\phi$ and 
$\psi$ generate the action of $\z^2$ on $W$.  We note that since $B$ 
is invariant for $\alpha$, $Y$ is invariant for $\psi$ as well as 
for $\phi$.

\proclaim{\thmee} $W/ (\id-\phi)W\cong Y/ (\id-\phi)Y$, and the 
isomorphism is equivariant for $\psi$.
\endproclaim

\demo{Proof} First we show that $W=(\id-\phi)W+Y$.  Let $x\in W$.  By 
\thmdd\ we have
$$\align
\phi(x)&=(\id-\phi)(x) + (2\phi-\id)(x) \\
&\in (\id-\phi)W+Y.
\endalign$$
Applying $\phi\inv$ we see that $x\in(\id-\phi)W+Y$.  Now we have
$${W\over(\id-\phi)W} = {(\id-\phi)W+Y\over(\id-\phi)W} \cong {Y\over 
Y\cap\bigl((\id-\phi)W\bigr)}.$$
We will show that $Y\cap\bigl((\id-\phi)W\bigr)=(\id-\phi)Y$, which 
will conclude the proof.  The containment ``$\supseteq$'' is clear.  
For the containment ``$\subseteq$'', let $y\in Y$ with 
$y=(\id-\phi)(x)$ for some $x\in W$.  Then
$\phi(x)=y+(2\phi-\id)(x)\in Y$, by \thmcc.  It follows that $x\in 
Y$, so that $y\in (\id-\phi)Y$. \qed
\enddemo

\proclaim{\thmff} $K_*(A)\cong K_*(B)$.
\endproclaim

\demo{Proof} \thmdd\ and \thmee, and the Pimsner-Voiculescu exact 
sequence, show that $K_*(A\times_\alpha\t^2\times_\phi\z)$ and 
$K_*(B\times_\alpha\t^2\times_\phi\z)$ are isomorphic, equivariantly for 
$\psi$.  Another application of Pimsner-Voiculescu, together with 
Takai-Takesaki duality, gives a commuting diagram of long exact sequences:
$$\minCDarrowwidth{.3 in}
\displaylines{
\CD
\cdots  @>>>
  {W\over(\id-\phi)W}  @>\id-\psi>>  {W\over(\id-\phi)W}  @>>> K_0(A)
  @>>>  \\
@.  {}  @AA{\cong}A  @AA{\cong}A  @AA{i_*}A    \\
\cdots  @>>>
  {Y\over(\id-\phi)Y}  @>\id-\psi>>  {Y\over(\id-\phi)Y}  @>>> K_0(B)
   @>>>  \\
\endCD\hfill\cr\cr\cr
\hskip1in
\CD
 @>>>  \ker\,(\id-\phi)  @>\id-\psi>>  \ker\,(\id-\phi)
   @>>> K_1(A) @>>> \cdots  \\
 @.  @|  @|  @AA{i_*}A  @.  \\
 @>>>  \ker\,(\id-\phi)  @>\id-\psi>>  \ker\,(\id-\phi)
   @>>> K_1(B) @>>> \cdots  \\
\endCD\hfill\cr\cr}
$$
It follows from the five lemma that $K_*(A)\cong K_*(B)$. \qed
\enddemo

\proclaim\thmkk Let $k\ge1$ be given.  For $0\le i$ and $1\le j\le k$
let $E_{i,j}$ be an irreducible directed graph with distinguished
vertex $v_{i,j}$ emitting infinitely many edges.  For $i\ge0$ let
$D_i$ be a copy of the graph $D$ in \thmj\ (with vertices $u_{i-1}$,
$u_i$, $a_i$, $a_i'$ --- see figure 2).  Let $\Omega$ be the
object obtained from the 1-graphs $\{D_i\}$ and the product $k$-graphs
$\{E_{i_1}\times\cdots\times E_{i,k}\}$ by identifying the vertex
$u_i$ with $(v_{i,1},\ldots,v_{i,k})$ as in \thmj.  Let $A$ be the
\cstar-algebra defined by the generators $\s$ and relations $\calr$ as
in \thmjj\ (modified in the obvious way).  Then $A$ is the unique UCT 
Kirchberg algebra with $K$-theory equal to
$$\bigoplus_{i=0}^\infty K_* \left(\otimes_{j=1}^k \oh(E_{i,j})\right).$$
\endproclaim

\demo{Proof} This follows from  \thmy\ and   \thmff.  (The 
uniqueness is a result of Zhang, \bib\zhang.) \qed
\enddemo

\bigskip
\epsfysize=2.5in                   
\centerline{\epsfbox{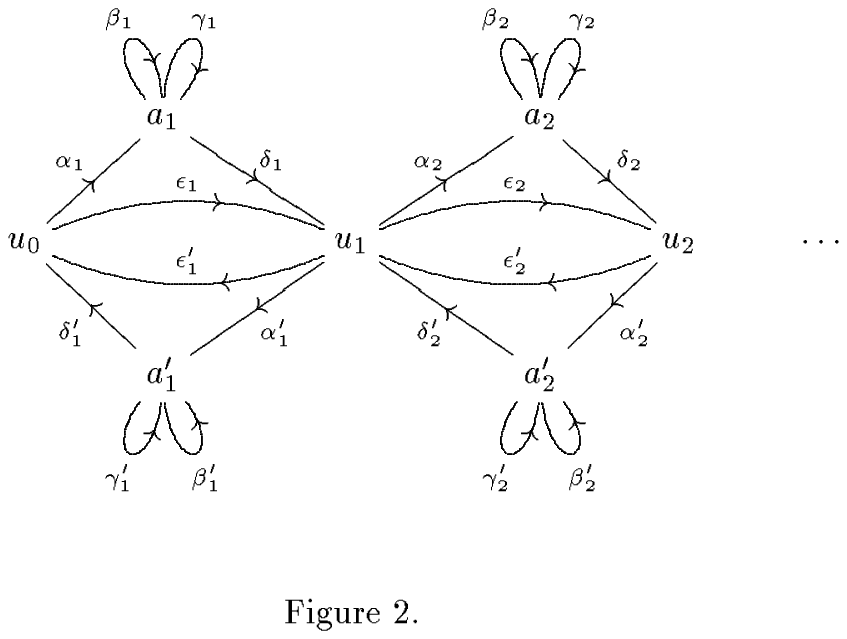}}
\bigskip

\head References \endhead

\roster
\item"\bib\anantharaman" C. Anantharaman-Delaroche, Purely infinite 
\cstar-algebras arising from dynamical systems, {\it Bull. Soc. Math. 
France\/} (1997), no. 125, 199-225.
\smallskip
\item"\bib\bkp" 
D.J. Benson, A. Kumjian and N.C. Phillips,
Symmetries of Kirchberg algebras, {\it Canad. Math. Bull.\/}
Vol. 46 (4) (2003), 509-528.
\smallskip
\item"\bib\cuntzon" J. Cuntz, Simple \cstar-algebras generated by
isometries, {\it Comm. Math. Phys.\/} {\bf 57} (1977), 173-185.
\smallskip
\item"\bib\cuntzoa" J. Cuntz, A class of
\cstar-algebras and topological Markov chains II: reducible chains and
the Ext-functor for \cstar-algebras, {\it Invent. Math.\/} {\bf 63}
(1981), 25-40.
\smallskip
\item"\bib\cuntzkrieger" J. Cuntz and W. Krieger, A class of
\cstar-algebras and topological Markov chains, {\it Invent. Math.\/}
{bf 56} (1980), 251-268.
\smallskip
\item"\bib\kirchberg" E. Kirchberg, The classification of purely 
infinite \cstar-algebras using Kasparov's theory, {\it Fields 
Institute Communications\/}, 2000.
\smallskip
\item"\bib \kumjianpask" A. Kumjian and D. Pask, Higher rank graph 
\cstar-algebras, {\it New York J. Math.\/} {\bf 6} (2000), 1-20.
\smallskip
\item"\bib\lacaspielberg" M. Laca and J. Spielberg,
Purely infinite \cstar-algebras 
from boundary actions of discrete groups, {\it J. reine angew. Math.\/} 
{\bf 480} (1996), 125-139.
\smallskip
\item"\bib\phillips" N.C. Phillips, A classification theorem for 
nuclear purely infinite simple \cstar-algebras, {\it Document math.\/} 
(2000), no. 5,  49-114.
\smallskip
\item"\bib\renault" J. Renault,  {\it A Groupoid Approach to 
\cstar-algebras\/},
Lecture Notes in Mathematics {\bf 793}, Springer-Verlag, Berlin, 1980.
\smallskip
\item"\bib\rordam" M. R\o rdam, Classification of Nuclear
\cstar-Algebras, in {\it Encyclopaedia of Mathematical Sciences\/}
{\bf 126} Springer-Verlag (2002), Berlin Heidelberg New York.
\smallskip
\item"\bib\spielberga" J. Spielberg, A functorial approach to the 
\cstar-algebras of a graph, {\it Internat. J. Math.\/} {\bf 13} 
(2002), no.3, 245-277.
\smallskip
\item"\bib\spielbergb" J. Spielberg, Semiprojectivity for certain 
purely infinite \cstar-algebras, pre\-print, Front for the
Mathematics ArXiv, math.OA/0102229.
\smallskip
\item"\bib\spielbergc" J. Spielberg, Non-cyclotomic 
presentations of modules and prime-order automorphisms of 
Kirchberg algebras, preprint, 2005.
\smallskip
\item"\bib\spielbergd" J. Spielberg, Weak semiprojectivity 
for purely 
infinite \cstar-algebras, pre\-print, Front for the
Mathematics ArXiv, math.OA/0206185 (revised 2005).
\smallskip
\item"\bib\szymanski" W. Szyma\'nski, The range of K-invariants for 
\cstar-algebras of infinite graphs,  {\it Indiana Univ. Math. J.}  
{\bf 51}  (2002),  no. 1, 239--249.
\smallskip
\item"\bib\zhang" S. Zhang, 
Certain \cstar-algebras with real rank zero and 
their corona and multiplier algebras, {\it Pacific J. Math.\/} {\bf 
155} (1992), 169-197.
\endroster

\enddocument
\bye